\documentclass[11pt]{article}
\usepackage{epsfig}
\usepackage{amssymb}
\usepackage{amscd}
\usepackage{epsf}
\usepackage{amsfonts,amssymb,amsthm,amsmath}

\long\def\comment#1\endcomment{\relax}

\newcounter{subsubsubsection}
\newcounter{subsubsubsubsection}

\makeatletter
\@addtoreset{subsubsubsection}{subsubsection}
\@addtoreset{subsubsubsubsection}{subsubsubsection}

\newcommand{\subsubsubsection}[1]{\par\addtocounter{subsubsubsection}{1}\smallskip
{\thesubsubsection.\arabic{subsubsubsection}.\bfseries\ #1}}

\makeatother

\newcommand{\sevafigc}[4]{\begin{figure}[h]\centerline{
 \epsfig{file=#1,width=#2,angle=#3}}
\bigskip\caption{#4}\end{figure}}

\DeclareMathSizes{11.1}{10}{8}{6}

\newtheorem*{theorem*}{Theorem}

\newtheorem*{conjecture}{Conjecture}

\newtheorem*{lemma}{Lemma}

\newtheorem*{klemma}{Key-Lemma}

\newtheorem*{proposition}{Proposition}

\newtheorem*{corollary}{Corollary}

\theoremstyle{remark}

\newtheorem*{remark}{Remark}

\newtheorem*{example}{Example}

\theoremstyle{definition}

\newcommand{\eqto}{\mathrel{\stackrel{\sim}{\to}}}

\newcommand{\mb}{{\bullet}}

\newcommand{\gl}{\mathfrak{gl}}

\def\wtilde#1{\widetilde{#1}\vphantom{#1}}

\title{ {\huge On $[A,A]/[A,[A,A]]$ and on a $W_n$-action \\on the
consecutive commutators of free \\associative algebra}}
\author{{\LARGE Boris Feigin and Boris Shoikhet}}

\begin{document}\maketitle

\begin{abstract}
\comment {\tt We consider the free associative algebra
$Free_n=<x_1,\dots,x_n>$ as a Lie algebra with the bracket
$[a,b]=a\cdot b-b\cdot a$. In this way we obtain a very complicated
infinite-dimensional Lie algebra $A_n$ which plays a crucial role in
the Assoc$\rightarrow$Lie version of the Kontsevich Formality. In
the present paper, we formulate a conjecture about the presentation
of $A_n$ in generators and relations. The generators are cyclic
expressions in $\{x_1,\dots, x_n\}$, and the relations are {\it not}
quadratic. This conjecture is checked on the MAGMA programm up to
the length of the words equal to 8 for $A_2$ and for ..... for
$A_3$.

We also consider in this paper the consecutive quotients $A_{n,k}$
of the lower central series of $A_n$, like
$A_{n,2}=[A_n,A_n]/[A_n,[A_n,A_n]]$ (here $k$ stands for the $k$-th
consecutive quotient). We consider $A_{n,k}$ as a $\gl_n$-module. It
is a graded $\gl_n$-module, $A_{n,k}=\oplus_{\ell\ge
1}A^{\ell}_{n,k}$ where $A^{\ell}_{n,k}$ is the subspace in
$A_{n,k}$ consisting of the monomials of length $\ell$. We
conjecture that for any $n$ and for any $k\ge 2$ the space
$A^{\ell}_{n,k}$ grows polynomially on $\ell$ for fixed $k$ for all
sufficiently large $\ell$. (The first quotient $A_n/[A_n,A_n]$ grows
exponentially on $\ell$). In particular, we prove that for any
$\ell$, $\dim A^{\ell}_{2,2}=\ell-1$, and $\dim
A^{\ell}_{3,2}=\ell^2-1$.}
\endcomment
We consider the lower central filtration of the free associative
algebra $A_n$ with $n$ generators as a Lie algebra. We consider the
associated graded Lie algebra. It is shown that this Lie algebra has
a huge center which belongs to the cyclic words, and on the quotient
Lie algebra by the center there acts the Lie algebra $W_n$ of
polynomial vector fields on $\mathbb{C}^n$. We compute the space
$[A_n,A_n]/[A_n,[A_n,A_n]]$ and show that it is isomorphic to the
space
$\Omega^2_{closed}(\mathbb{C}^n)\oplus\Omega^4_{closed}(\mathbb{C}^n)\oplus\Omega^6_{closed}(\mathbb{C}^n)
\oplus\dots$.
\end{abstract}
\comment The {\tt Assoc}$\rightarrow${\tt Lie} version of the
Kontsevich Formality (also conjectured by Maxim Kontsevich) deals
with the "second fundamental functor" which associates with an
associative algebra $A$ the Lie algebra {\tt Lie}$(A)$ which is $A$
as a vector space, and with the bracket $[a,b]=a\cdot b-b\cdot a$
(where $\cdot$ is the product in $A$). Consider the case where
$A=Free_n=<x_1,\dots,x_n>$, and put $A_n={\tt Lie}(A)$. Then the
{\tt Assoc}$\rightarrow${\tt Lie} formality conjecture states that
the Lie cochain complex with the adjoint coefficients
$C^\mb(A_n,A_n)$ is formal as differential graded Lie algebra with
the Gerstenhaber bracket. This complex plays the role of the
Hochschild cohomological complex in the case of Kontsevich formality
theorem.

The main difficulty in the {\tt Assoc}$\rightarrow${\tt Lie}
formality is to compute the Lie algebra cohomology
$H_{Lie}(A_n,A_n)$, that is, the analog of polyvector fields. In the
present paper we conjecturally describe the Lie algebra $A_n$ in
generators and relations. This can be considered as the first step
in the computation of this cohomology.
\endcomment

\section*{Introduction}
Let $A$ be an associative algebra. A free resolution
$\mathcal{R}^\mb$ of $A$ is a free graded differential algebra
$\mathcal{R}^\mb=\bigoplus R^i$, $i\in\mathbb{Z}_{\le 0}$, with
differential $Q$ has degree $+1$ such that the cohomology of $Q$ is
only in degree zero, and is canonically isomorphic to $A$ as
algebra. Such a resolution can be used for calculation of "higher
derived functors" for $A$. For example, higher cyclic homology of
$A$ is the higher derived functor for the functor
$$A\to A/[A,A].$$
It means that for the calculation of cyclic homology of $A$ we have
to take an arbitrary free resolution $\mathcal{R}^\mb$ of $A$ and
consider the quotient
$\mathcal{R}^\mb/[\mathcal{R}^\mb,\mathcal{R}^\mb]$. The
differential $Q$ acts in
$\mathcal{R}^\mb/[\mathcal{R}^\mb,\mathcal{R}^\mb]$, and cohomology
of $Q$ is the higher derived functor of  $A\to A/[A,A].$

It is natural to try to calculate "higher derived" for other
functors. Surely there is a lot of interesting and important
functors, but unfortunately in the most cases it is rather hard to
calculate the higher derived functors. The cyclic homology is an
exception, because it can be expressed in the terms of usual
Hochschild homology. Another relatively simple case is the functor
of abelianization
$$A\to A^{ab}\simeq A/J,$$
where $J$ is the two-sided ideal generated by the brackets $[a,b]$,
$a,b\in A$. Let $A$ be the polynomial algebra
$\mathbb{C}[x_1,\ldots,x_n]$. A resolution $\mathcal{R}^\mb$ of $A$
can be constructed in terms of the dual Grassmannian algebra
$\Lambda^\mb(\xi_1,\ldots,\xi_n)$. Let $A_+^\vee$ be the kernel of
the augmentation map
$$\Lambda^\mb(\xi_1,\ldots,\xi_n)\to\mathbb{C}.$$
Then $A_+^\vee$ is a graded algebra and the resolution
$\mathcal{R}^\mb$ of $A$ is a free graded algebra generated by the
dual shifted space $(A_+^\vee)^*[-1]$. The differential in
$\mathcal{R}^\mb$ is given by the coproduct
$$(A_+^\vee)^*\to (A_+^\vee)^*\otimes (A_+^\vee)^*$$.
On the quotient
$$\mathcal{R}^\mb/\mathcal{R}^\mb[\mathcal{R}^\mb,\mathcal{R}^\mb]\mathcal{R}^\mb$$ the
differential acts by zero. Therefore the "higher abelianization" of
$\mathbb{C}[x_1,\ldots,x_n]$ is an algebra of functions on super
vector space $A_+^\vee[-1]$ (for example, the "higher
abelianization" of $\mathbb{C}[x_1,x_2]$ is the algebra
$\mathbb{C}[u_1,u_2,u_{1,2}]$, $\deg u_i=0$, $\deg u_{1,2}=-1$).

In this paper we study the functor
$$A\to A/[A,[A,A]].$$
In order to say something about the higher functors we need to know
what is the quotient  $A/[A,[A,A]]$ for a free algebra $A$. The
related question is to determine the higher functors for
$$
A\to A/A[A,[A,A]]A
$$

Our first result is an explicit calculation of the last quotient for
the free algebra $T(V)$ of an $n$-dimensional vector space $V$; we
also denote this algebra by $A_n$. Let
$$\Omega^\mb=S^\mb(V)\otimes \Lambda^\mb(V)$$
be the de Rham complex of the polynomial differential forms on the
dual space $V^*$. Differential $d$ on $\Omega^\mb$ determines the
bivector field $d\wedge d=\nu$. Note that
$$d^2=\frac12[d,d]=0,$$
therefore $[\nu,\nu]=0$ and $\Omega^\mb$ with the bracket
$$[\omega_1,\omega_2]=(-1)^{\deg \omega_1}\cdot 2 d\omega_1 \wedge d\omega_2$$
is a $\mathbb{Z}_2$-graded Poisson algebra. A quantization of
$(\Omega,\nu)$ may be given  by the very simple formula:
$$\omega_1*\omega_2=\omega_1\wedge\omega_2+ (-1)^{\deg \omega_1}
d\omega_1 \wedge d\omega_2.$$ Denote the result of quantization by
$\Omega(V)_*$. Then $\Omega(V)_*$ is a $\mathbb Z_2$-graded
associative algebra and the algebra of even forms with the quantized
product $\Omega^{even}(V)_*$ is a subalgebra of $\Omega(V)_*$
generated by the space $V\hookrightarrow \Omega^0(V)$.

\begin{proposition}
The quotient algebra
$$B=A_n/A_n[A_n,[A_n, A_n]]A_n$$
is isomorphic to $\Omega^{even}(V)_*$. The map
$B\to\Omega^{even}(V)_*$ restricted to the space of generators $V$
is the identity map.
\end{proposition}

The second commutator
$$[\Omega^{even}(V)_*,[\Omega^{even}(V)_*, \Omega^{even}(V)_*]]$$
vanishes, so the natural homomorphism $A_n\to \Omega^{even}(V)_*$
induces a map
$$\theta:[A_n,A_n]/[A_n,[A_n,A_n]]\to [\Omega^{even}(V)_*,\Omega^{even}(V)_*].$$
It is easy to see that  $[\Omega^{even}(V)_*,\Omega^{even}(V)_*]$
coincides with the space of exact (=closed of degree $>0$) even
forms. We prove the map $\theta$ is an isomorphism
$$
\theta: [A_n,A_n]/[A_n,[A_n,A_n]]\to \Omega^{even>0}_{closed}(V)
$$
This result is equivalent to the fact that
$$
[A_n,[A_n,A_n]]=[A_n,A_n]\cap\left(A_n\cdot [A_n,[A_n,A_n]]\cdot
A_n\right)
$$

Let us summarize. We consider the following functors on the category
of associative  algebras:
$$
\begin{aligned}
\ &F_1: B\to B/[B,B],\\
 &F_2: B\to B/[B,[B,B]],\\
&G_1: B\to B/B[B,B]B,\\
&G_2: B\to B/B[B,[B,B]]B,\\
&F_{1,2}: B\to [B,B]/[[B,B],B],\\
&G_{1,2}: B\to B/(B[B,[B,B]]B + [B,B])
\end{aligned}
$$
Higher derived functors for $F_1, F_2, G_1, G_2, F_{1,2}, G_{1,2}$
can be found without a lot of problems when we know the higher
abelianization and cyclic homology.

For example, let $B=\mathbb{C}[x_1,x_2]$. We know that a free
resolution of $\mathbb{C}[x_1,x_2]$ is the algebra generated by
$u_1,u_2$ and $u_{1,2}$. Differential $Q$ is defined by
$$Q(u_i)=0,\ Q(u_{1,2})=u_1u_2-u_2u_1.$$
Functor $G_2$ applied to the resolution gives us the even part of
the algebra of forms as a superspace:
$$\mathbb{C}[u_1,u_2,u_{1,2}; du_1,du_2,du_{1,2}].$$ The
differential $Q$ acts nontrivially only on $u_{1,2}$:
$Qu_{1,2}=du_1\wedge du_2$. So the higher derived functor for the
functor $G_2$ is just the cohomology of this differential.

Now consider the exact sequence of functors
$$0\to F_{1,2}\to F_2\to F_1\to 0.$$
Using it we can find "higher derived" for $F_2$.

In the beginning of this work we tried to analyze the numerical
results on the dimensions $A_{n,k}^\ell$ of graded components for
algebra $A_n$ for $n=2,3$ and small $k,\ell$, obtained by Eric Rains
on MAGMA (see (17) below). Here $A_{n,k}$ is the quotient space of
$k$-commutators of $A_n$ modulo $k+1$-commutators, and
$A_{n,k}^\ell$ is the component of $A_{n,k}$ of monomials of the
length $\ell$. For example, $A_{n,1}=A_n/[A_n,A_n]$ is the space of
cyclic words on $n$ variables, and the dimensions $A_{n,1}^\ell$
grow exponentially as $\ell$ tends to $\infty$. Our first
observation was an unexpected phenomena that the spaces
$A_{n,k}^\ell$ for $k\ge 2$ grow {\it polynomially} on $\ell$. We
saw it from (17) for small $k,\ell$ and $n=2,3$. In general, it is a
conjecture till now.

Then, if they grow polynomially, we tried to think about them as
about some tensor fields, more precisely, (maybe not irreducible)
$W_n$-modules. The present paper is the result of our attempt to
understand these two phenomena--the polynomial growth of
$A_{n.k}^\ell$, and a structure of $W_n$-module on it.

Another strange thing which appeared from our results is that the
space $[A_n,A_n]/[A_n,[A_n,A_n]]$ is an {\it associative algebra}.
The algebra structure is very unclear from this definition, but it
follows from the description of the last space as
$\Omega^{even>0}_{closed*}(V)$. It is a commutative algebra with the
usual wedge product of differential forms.

At the moment we can not find such a theory for higher $A_{n,k}$,
$k>2$. We have some conjectures which hopefully will be published
somewhere.

Now let us outline the contents of the paper:

In Section 1 we prove "by hands" that $[A_n,A_n]/[A_n,[A_n,A_n]]$ is
isomorphic to the space of closed 2-forms on $\mathbb{C}^n$ for
$n=2,3$;

In Section 2 we develop our main technics, aimed to find
$A_n/(A_n[A_n,[A_n,A_n]]A_n)$ and $[A_n,A_n]/[A_n,[A_n,A_n]]$ for
general $n$, we prove our results here modulo Lemma 2.2.2.2 which is
proven in Section 3. To prove this Lemma in Section 3 we use a
theorem describing all irreducible $W_n$-modules of a reasonable
class;

In Section 4 we define a $W_n$-action on the quotient (by the
center) $\wtilde{gr}(A_n)$ of the associated graded Lie algebra of
$A_n$ with respect to the lower central filtration. We also consider
many examples here. \comment
\subsection*{Acknowledgements.}We are grateful to Pavel Etingof, Giovanni
Felder, Anton Khoroshkin,  and especially to Misha Movshev and to
Eric Rains for many discussions. B.F. is grateful to the grants RFBR
05-01-01007, SS 2044.2003.2, RFBR-JSPS 05-01-02934, and also to the
Russian Academy of Science project "Elementary particle physics and
fundamental nuclear physics" for partial financial support. B.Sh. is
thankful to the ETH (Zurich), and to the University of Luxembourg,
where some ideas of this paper were invented, for excellent working
conditions. He is also grateful to the ETH Research Comission grant
TH-6-03-01 and to the research grant R1F105L15 of the University of
Luxembourg for partial financial support.
\endcomment

\comment
\section{Generators and relations}
\subsection{The Hilbert polynomials $H_n$ and $H_n^\prime$}
Let $A_n$ be the free associative algebra with generators
$\{x_1,\dots, x_n\}$, considered as a Lie algebra with the bracket
$[a,b]=a*b-b*a$. For any pronilpotent Lie algebra $A$ (and, in
particular, for $A_n$) the space of generators is $A/[A,A]$. In the
case of the Lie algebra $A_n$ this quotient is isomorphic to the
space of {\it cyclic words} in $\{x_1,\dots,x_n\}$.

Consider the lower central series of $A_n$: $F_1=A_n$, and
$F_{k+1}=[A_n,F_k]$. Denote $A_{n,k}=A_k=F_k/F_{k+1}$, and consider
the associated graded Lie algebra {\tt gr}$A_n=\oplus_{k\ge
1}A_{n,k}=\oplus_{k\ge 1}F_k/F_{k+1}$. Then the space
$A_n/[A_n,A_n]$ is {\it canonically} imbedded to {\tt gr}$A_n$.

There is also the canonical space of quadratic "relations" which is
the kernel of the map $\Lambda^2(F_1/F_2)\to F_2/F_3=A_{n,2}\subset
{\tt gr}A_n$. Namely, suppose $A$ is an associative algebra, and
$A={\tt Lie}(A)$. Suppose $A/[A,A]$ is the space of generators. For
$a\in A$ denote by $\overline{a}$ its image in $A/[A,A]$. Then the
following relations hold:
\begin{equation}\label{eq1_1}
[\overline{a_1},\overline{a_2\cdot
a_3}]+[\overline{a_2},\overline{a_3\cdot
a_1}]+[\overline{a_3},\overline{a_1\cdot a_2}]=0
\end{equation}
in $[A,A]/[A,[A,A]]$ for any $a_1,a_2,a_3\in A$, because there holds
$$
[a_1,a_2\cdot a_3]+[a_2,a_3\cdot a_1]+[a_3, a_1\cdot a_2]=0
$$
in $A$.

We can {\it not} consider the relations (\ref{eq1_1}) as relations
between generators in $A$, but only in ${\tt gr}A$. Following this
line, we firstly studied the following two questions:

\begin{itemize}
\item[1.] Is the Lie algebra ${\tt gr}A_n$ isomorphic to
$A_n$?
\item[2.] Is the Lie algebra ${\tt gr}A_n$ quadratic
with relations (\ref{eq1_1})?
\end{itemize}
\begin{example}{\tt The case $n=1$.}
\end{example}
Consider the case $n=1$. In this case the Lie algebra $A_1$ is
Abelian and is, therefore, isomorphic to ${\tt gr}A_1$. Let us show
that this Lie algebra is quadratic with relations (\ref{eq1_1}). The
monomials $\{1,x,x^2,\dots\}$ form a basis in $A_1$ as in a vector
space. The relations (\ref{eq1_1}) can be rewritten as
\begin{equation}\label{eq1_1new}
[x^k,x^{N-k}]+[x^\ell, x^{N-\ell}]+[x^m,x^{N-m}]=0
\end{equation}
for any $N\ge 0$, $k+\ell+m=N$. Consider the free Lie algebra
generated by $\{1,x,x^2,\dots\}$ and with the relations
(\ref{eq1_1new}). Then the system of linear equations on $[x^a,x^b]$
is overdefined, and it implies that each bracket $[x^a,x^b]$ is 0.

\bigskip
 We show now by some computations that the answer for $n>1$ for the
second question is negative. (The answer for the first question is
probably also negative, but it does not follow from our
computations). The computations below are made using MAGMA.

Consider the bigraded Hilbert series for $A_n$:
\begin{equation}\label{eq1_2}
H_{n}=\sum_{\ell\ge 0,\ k\ge 1}
\dim {A_{n,k}^\ell}u^k t^\ell
\end{equation}
For $n=2$ and for $n=3$ the bigraded Hilbert series are:
\begin{equation}\label{eq1_3}
\begin{aligned}
\ & H_2(u,t)=\\
&\ \ \ (  u)\\
& +( 2u)t\\
& +( 3u+ u^2)t^2\\
& +( 4u+2u^2+ 2u^3)t^3\\
& +( 6u+3u^2+ 4u^3+ 3u^4)t^4\\
& +( 8u+4u^2+ 6u^3+ 8u^4+ 6u^5)t^5\\
& +(14u+5u^2+ 8u^3+13u^4+15u^5+ 9u^6)t^6\\
&+(20u+6u^2+10u^3+18u^4+26u^5+30u^6+ 18u^7)t^7\\
&+(36u+7u^2+12u^3+23u^4+37u^5+54u^6+ 57u^7+ 30u^8)t^8\\
&+(60u+8u^2+14u^3+28u^4+48u^5+80u^6+108u^7+110u^8+56u^9)t^9\\
&+\mathcal{O}(t^{10})\\
& \\
&\\
 &H_3(u,t)=\\
 &\ \ \  (   u)\\
&+(  3u)t\\
& +(  6u+ 3u^2)t^2\\
& +( 11u+ 8u^2+ 8u^3)t^3\\
& +( 24u+15u^2+24u^3+ 18u^4)t^4\\
& +( 51u+24u^2+48u^3+ 72u^4+ 48u^5)t^5\\
&+(130u+35u^2+80u^3+162u^4+206u^5+116u^6)t^6\\
&+\mathcal{O}(t^7)
\end{aligned}
\end{equation}
\begin{remark}
We see from these formulas that the coefficients in $u^2$ and $u^3$
depend polynomially on the exponent of $t$. The same is true for
$u^4$,... but starting from some place. Comparably, the coefficients
in $u^1$ depend exponentially in the exponent of $t$. We will
discuss these things in the next Section.
\end{remark}

Consider now the following quadratic Lie algebra $A^{\prime}_n$: it
is generated by $A_n/[A_n,A_n]$ with the quadratic relations
(\ref{eq1_1}). Clearly there is a map of graded Lie algebras
$\varphi\colonA_n^\prime\to{\tt gr}A_n$. It is surjective by the
definition. We can associate with $A_n^\prime$ the bigraded Hilbert
series
\begin{equation}\label{eq1_4}
H^{\prime}_{n}=\sum_{\ell\ge 0,\ k\ge 1} \dim
{(A_{n,k}^\ell)^\prime} u^k t^\ell
\end{equation}
where $(A_{n,k}^\ell)^\prime$ is the graded component which has the
grading $k$ by generators and the grading $\ell$ by the length of
monomial in the alphabet $\{x_1,\dots,x_k\}$. Our goal is to compare
the Hilbert polynomials $H_n$ and $H_n^\prime$. We will see that
$H_n^\prime$ is strictly greater than $H_n$ even for 2 generators
(for 1 generator they are the same--see Example above).

For $n=2$ and $n=3$ the bigraded Hilbert series $H_n^\prime$ are:

\begin{equation}\label{eq1_5}
\begin{aligned}
\ & H^{\prime}_2(u,t)=\\
& \ \ \ (  u)\\
& +( 2u)t\\
& +( 3u+ u^2)t^2\\
& +( 4u+2u^2+ 2u^3)t^3\\
& +( 6u+3u^2+ 4u^3+ 3u^4)t^4\\
& +( 8u+4u^2+ 6u^3+ 8u^4+ 6u^5)t^5\\
& +(14u+5u^2+ 8u^3+13u^4+16u^5+ 9u^6)t^6\\
&+(20u+6u^2+10u^3+18u^4+32u^5+ 32u^6+ 18u^7)t^7\\
&+(36u+7u^2+12u^3+23u^4+48u^5+ 70u^6+ 64u^7+ 30u^8)t^8\\
&+(60u+8u^2+14u^3+28u^4+64u^5+122u^6+160u^7+128u^8+56u^9)t^9\\
&+\mathcal{O}(t^{10})\\
& \\
& \\
& H^{\prime}_3(u,t)=\\
&\ \ \  (   u)\\
&+(  3u)t\\
& +(  6u+ 3u^2)t^2\\
& +( 11u+ 8u^2+ 8u^3)t^3\\
& +( 24u+15u^2+24u^3+ 18u^4)t^4\\
& +( 51u+24u^2+48u^3+ 72u^4+ 48u^5)t^5\\
&+(130u+35u^2+80u^3+172u^4+216u^5+116u^6)t^6\\
&+\mathcal{O}(t^7)
\end{aligned}
\end{equation}

The differences $H_n^\prime-H_n$ for $n=2$ and $n=3$ are
\begin{equation}\label{eq1_6}
\begin{aligned}
\ & H_2^\prime(u,t)-H_2(u,t)=(  u^5)t^6 +( 6u^5+ 2u^6)t^7+
(11u^5+16u^6+ 7u^7)t^8\\
&+(16u^5+42u^6+52u^7+18u^8)t^9+\mathcal{O}(t^{10}),\\
&\\
&H_3^\prime(u,t)-H_3(u,t)=(10u^4+10u^5)t^6+\mathcal{O}(t^7)
\end{aligned}
\end{equation}

For 4 and 5 generators we can not see the difference because the
computer can compute only up to $t^4$.

We see also that it is impossible to see "by hands" that
$H_n^\prime-H_n$ is not zero: for $n=2$ the first difference appears
when the length of the monomial is 6, and it belongs to
$F_5^\prime/F^\prime_6$!
\subsection{The Conjecture}
We see from the computation above that the relations (\ref{eq1_1})
are {\it all} quadratic relations on the generators, that is
\subsubsection{}
\begin{lemma}
\begin{equation}\label{eq1_7}
\left(\Lambda^2(A_n/[A_n,A_n])\right)/({\tt relations\
(1)})=[A_n,A_n]/[A_n,[A_n,A_n]]
\end{equation}
\begin{proof}
The following proof is almost contained in the Loday's book on
Cyclic Homology [L].

Let $A$ be an associative algebra, and $A={\tt Lie}(A)$. Denote
$V(A)=\left(\Lambda^2(A)\right)/(a\wedge (b\cdot c)+b\wedge (c\cdot
a)+c\wedge(a\cdot b)=0)$ for all $a,b,c\in A$. We have the following
short exact sequence of $A$-modules (with respect to the adjoint
action):
$$
0\longrightarrow HC_1(A)\longrightarrow V(A)\longrightarrow
[A,A]\longrightarrow 0
$$
where $HC_1(A)$ is the first cyclic homology of the associative
algebra $A$.

Consider the corresponding long exact sequence in Lie algebra
homology:
$$
\dots\rightarrow HC_1(A)\rightarrow
\left(\Lambda^2(A/[A,A])\right)/(a\wedge (b\cdot c)+b\wedge (c\cdot
a)+c\wedge(a\cdot b)=0)\rightarrow [A,A]/[A,[A,A]]\rightarrow 0
$$
One should explain the middle term. It is not true that
$\Lambda^2(A)/[A,\wedge^2(A)]$ is $\Lambda^2(A/[A,A]$, and we should
use the identity $a\wedge (b\cdot c)+b\wedge (c\cdot
a)+c\wedge(a\cdot b)=0$. We have by this identity:
$$
[a,b\wedge c]=[a,b]\wedge c+b\wedge [a,c]=a\wedge [b,c]
$$
which explains the middle term.

Now apply the last long exact sequence to the case when
$A=Free<x_1,\dots, x_n>$. Then $HC_1(A)=0$, and we obtain the
isomorphism (\ref{eq1_7}).
\end{proof}
\end{lemma}
\bigskip
Then the situation is as follows: the Lie algebra ${\tt gr}A_n$ is
not quadratic, and one can expect relations of higher degree. (We
know by the argument above that the relations (\ref{eq1_1}) exhaust
all quadratic relations). The problem of finding these higher
relations in ${\tt gr}A_n$ seems to be very complicated.

\bigskip
\subsubsection{} On the other hand, we now modify the relations
(\ref{eq1_1}) to get the relations in the initial Lie algebra $A_n$.

We should firstly imbed the abstract cyclic words $A_n/[A_n,A_n]$ to
$A_n$. We choose the following imbedding: we map a monomial
$m=x_{i_1}x_{i_2}\dots x_{i_k}$ to the cyclic sum
\begin{equation}\label{eq1_8}
\wtilde{m}=x_{i_1}x_{i_2}\dots x_{i_k}+x_{i_2}\dots
x_{i_k}x_{i_1}+\dots +x_{i_k}x_{i_1}x_{i_2}\dots
\end{equation}
\begin{example}
For $m=x_1^2x_2x_3$,
$\wtilde{m}=x_1^2x_2x_3+x_1x_2x_3x_1+x_2x_3x_1^2+x_3x_1^2x_2$
\end{example}
\begin{lemma}
For any monomial $m\in A_n$ of the length $\ell$ one has
\begin{equation}\label{eq1_9}
m=\frac1\ell\wtilde{m}+[A_n,A_n]
\end{equation}
\begin{proof} It is clear.
\end{proof}
\end{lemma}
Consider now equation (\ref{eq1_1}) as follows:
\begin{equation}\label{eq1_10}
[a,b\cdot c]+[b,c\cdot a]+[c, a\cdot b]=0
\end{equation}
where $a,b,c$ are considered as elements of an associative algebra
$A$, and {\it not} as elements of $A/[A,A]$ , as we considered in
Section 1.1. Then (\ref{eq1_10}) is true for any associative algebra
$A$ and for the Lie bracket on $A={\tt Lie}(A)$.

Our goal now is to replace in (\ref{eq1_10}) $a$ for $\wtilde{a}$,
$b$ for $\wtilde{b}$, and $c$ for $\wtilde{c}$. For this we use
Lemma above.

We obtain from (\ref{eq1_10}) equation of the form
\begin{equation}\label{eq1_11}
\begin{aligned}
\ &\ \ \
\left[\frac{1}{\ell(a)}\wtilde{a}+\sum[a_1^i,a_2^i],\frac{1}{\ell(b\cdot
c)}\wtilde{b\cdot
c}+\sum[\alpha_1^j,\alpha_2^j]\right]\\
&+\left[\frac{1}{\ell(b)}\wtilde{b}+\sum[b_1^i,b_2^i],\frac{1}{\ell(c\cdot
a)}\wtilde{c\cdot
a}+\sum[\beta_1^j,\beta_2^j]\right]\\
&+\left[\frac{1}{\ell(c)}\wtilde{c}+\sum[c_1^i,c_2^i],\frac{1}{\ell(a\cdot
b)}\wtilde{a\cdot b}+\sum[\gamma_1^j,\gamma_2^j]\right]=0
\end{aligned}
\end{equation}
Now we apply Lemma above to all $a_1^i,a_2^i,b_1^i,b_2^i,
c_1^i,c_2^i$ and to
$\alpha_1^j,\alpha_2^j,\beta_1^j,\beta_2^j,\gamma_1^j,\gamma_2^j$,
etc. Finally, we obtain a relation which contains only cyclic
expressions, but highly non-quadratic in brackets. In this way we
obtain a cyclic relation $R(a,b,c)$ depending on a triple of
elements of $A$. Now we are ready to formulate our main Conjecture:
\begin{conjecture}
The Lie algebra $A_n$ is generated by the cyclic words of the form
$\wtilde{a}$ for the monomials $a\inA_n$, and with relations
$R(a,b,c)$ for each triple of monomials $a,b,c\inA_n$.
\end{conjecture}

\bigskip
Philosophically, in the relations $R(a,b,c)$ we may have the same
quadratic part for $(a,b,c)$ and $(a^\prime,b^\prime, c^\prime)$,
and then the difference $R(a,b,c)-R(a^\prime, b^\prime, c^\prime)$
is a cubic relation, etc.

The conjecture is checked up to the length of monomials equal to 8
for $n=2$ and equal to .... for $n=3$.
\endcomment

\section{The isomorphism $[A_n,A_n]/[A_n,[A_n,A_n]]\simeq
\Omega^2_{closed}(\mathbb{C}^n)$ for \\$n=2,3$} In this Section we
compute "by hands" the quotient $[A_n,A_n]/[A_n,[A_n,A_n]]$ for
$n=2$ and $n=3$. To make the exposition more clear, we first define
the concept of a non-commutative 1-form.
\subsection{Non-commutative 1-forms}
Let $A$ be an associative algebra. A 1-form on $A$ is a finite sum
of the expressions $a\cdot db\cdot c$, where $a,b,c\in A$ modulo the
following two relations:
\begin{itemize}
\item[(i)] the cyclicity: $t\cdot a\cdot db\cdot c=a\cdot db\cdot
c\cdot t$ for any $a,b,c,t\in A$,
\item[(2)] the Leibniz rule: $d(a\cdot b)=(da)\cdot b+a\cdot db$.
\end{itemize}
We say that a 1-form on $A$ is exact if it has the form $\omega=da$,
$a\in A$.

Consider the space $\Omega^1_A/d\Omega^0_A$ of 1-forms on $A$ modulo
the exact 1-forms. We can reduce any 1-form to an expression $a\cdot
db$ using the cyclicity. Next, modulo the exact forms $a\cdot
db+b\cdot da=0$. Thus, there is a map of the space
$\Lambda^2(A)\to\Omega^1_A/d\Omega^0_A$ which is clearly surjective.
What is its kernel?

We have the relation:
$a_1d(a_2a_3)=a_1d(a_2)a_3+a_1a_2d(a_3)=a_3a_1d(a_2)+a_1a_2d(a_3)$
which is
\begin{equation}\label{eqrev8}
a_1\wedge(a_2a_3)=(a_3a_1)\wedge a_2+(a_1a_2)\wedge a_3
\end{equation}
or, in more symmetric form,
\begin{equation}\label{eqrev9}
(a_1a_2)\wedge a_3+(a_2a_3)\wedge a_1+(a_3a_1)\wedge a_2=0
\end{equation}
It is clear that there are no other relations. We proved the
following result:
\begin{lemma}
For any associative algebra $A$, the space $\Omega^1_A/d\Omega^0_A$
is isomorphic to $\Lambda^2(A)/({\it relations} (2))$.
\end{lemma}
\qed \comment Consider an associative algebra $A$ with unit. In [K2]
M.Kontsevich defines the functions associated with $A$ as $A/[A,A]$.
Notice that it is not an algebra, only a vector space. Next, he
defines the de Rham complex associated with $A$ as follows:

Firstly, consider the differential envelope $DA$ of the algebra $A$.
This is by definition the algebra generated by the elements $a$ and
$da$, where $a\in A$, and we have: $d(ab)=(da)b+adb$. A general
element of $DA$ has the form $\omega=a_1(da_2)a_3(da_4)a_5\dots$.
The algebra $DA$ is a differential $\mathbb{Z}$-graded algebra. The
differential maps $a$ to $da$ and $da$ to $0$, and satisfies the
Leibniz rule. We call it the de Rham differential.  Finally, define
the differential forms associated with $A$ as
$\Omega^\mb(A)=DA/[DA,DA]$. Again, it is not an algebra, only a
$\mathbb{Z}$-graded vector space. We still have the de Rham
differential on $\Omega^\mb(A)$, because the de Rham differential
preserves $[DA,DA]$.

Consider the case where $A=A_n$, the free associative algebra with
$n$ generators. Then we have the Poincare lemma: $\Omega^\mb$ is
acyclic except degree 0 where the cohomology is 1-dimensional. To
prove it, first define a vector field on an $n$-dimensional
non-commutative space as a derivation of the algebra $A_n$. Let
$\xi$ be such a vector field. Then there are defined two "vector
fields" acting on $\Omega^\mb(A_n)$. These are the Lie derivation
$L_\xi$ and the insertion operator $i_\xi$, defined on $DA$ as
$L_\xi(a)=\xi(a)$, $L_\xi(da)=d(\xi(a))$, and $i_\xi(a)=0$,
$i_\xi(da)=\xi(a)$ for any $a\in A$, and then we extend them to the
action on $DA$ by the Leibniz rule. The both operators $L_\xi$ and
$i_\xi$ preserve $[DA,DA]$, and therefore act on $\Omega^\mb(A)$.
They satisfy the following relations: $L_\xi=i_\xi d+di_\xi$, $i_\xi
i_\eta+i_\eta i_\xi=0$, $[L_\eta,i_\xi]=i_{[\eta,\xi]}$,
$L_{[\xi,\eta]}=[L_\xi, L_\eta]$. Consider the vector field
$e=x_1\frac{\partial}{\partial
x_1}+\dots+x_n\frac{\partial}{\partial x_n}$ on $A_n$. We have:
$L_e=i_ed+di_e$, and $L_e$ is the length-grading operator. This
gives us the desired homotopy.

In particular, for $A=A_n$, the closed forms are the same that the
exact forms. We denote the de Rham complex for $A=A_n$ by
$\Omega^\mb_{nc}$ omitting $n$.

Now we describe here $\Omega_{nc,\ closed}^2=d\Omega^1_{nc}$.

Notice that this space is isomorphic to the space
$\Omega^1_{nc}/d\Omega^0_{nc}$ by the Poincare lemma. We describe
the latter space.

Any $1-form$ can be written as $adb$, where $a,b\in A$. In the
quotient space we have: $d(ab)=0$. But $d(ab)=(da)b+adb$. The first
summand is equal to $bda$ because we define differential forms as
$DA/[DA,DA]$. In this way we obtain $\Lambda^2(A)$. We have one more
relation:
$a_1d(a_2a_3)=a_1d(a_2)a_3+a_1a_2d(a_3)=a_3a_1d(a_2)+a_1a_2d(a_3)$.
We obtain the following relation:
\begin{equation}\label{eqrev8}
a_1\wedge(a_2a_3)=(a_3a_1)\wedge a_2+(a_1a_2)\wedge a_3
\end{equation}
or, in more symmetric form,
\begin{equation}\label{eqrev9}
(a_1a_2)\wedge a_3+(a_2a_3)\wedge a_1+(a_3a_1)\wedge a_2=0
\end{equation}
It is clear that these are the only relations for
$\Omega^1_{nc}/d\Omega^0_{nc}$. We prove the following fact:

The space $\Omega^1_{nc}/d\Omega^0_{nc}$ is isomorphic to
$\left(\Lambda^2(A_n)\right)/{\tt ({\it relations}\
(\ref{eqrev9}))}$.

Now we can prove the following result:
\begin{lemma}{\rm(Kontsevich [K2])}
The space $\Omega^2_{nc,\ closed}$ is isomorphic to the space
$[A_n,A_n]$.
\begin{proof}
We have the commutator map
$[,]\colon\left(\Lambda^2(A_n)\right)/{\tt ({\it relations}\
(\ref{eqrev9}))}\to[A_n,A_n]$. It is a part of the short exact
sequence
$$
0\longrightarrow HC_1(A)\longrightarrow \left(\Lambda^2
A\right)/{\tt ({\it relations}\  (\ref{eqrev9}))}\longrightarrow
[A,A]\longrightarrow 0
$$
which holds for any associative algebra $A$. Now the result follows
from the fact that the first cyclic homology $HC_1(A_n)$ is 0.
\end{proof}
\end{lemma}
\begin{remark}
In [K2], Kontsevich claims a different proof, using a computation of
Euler characteristics.
\end{remark}
\endcomment
\subsection{A Lemma}
In the case when $A=A_n$, the free associative algebra with $n$
generators over $\mathbb{C}$, the space $\Lambda^2(A_n)/({\it
relations} (2))$ is isomorphic to the commutator $[A_n,A_n]$. We
have the following Lemma:
\begin{lemma}
\begin{itemize}
\item[(i)]
$\Lambda^2(A_n)/({\it relations} (2))=[A_n,A_n]$,
\item[(ii)]
$\left(\Lambda^2(A_n/[A_n,A_n])\right)/({\tt {\it relations}\
(2)})=[A_n,A_n]/[A_n,[A_n,A_n]]$
\end{itemize}
\begin{proof}
For any associative algebra $A$, we have the short exact sequence:
\begin{equation}\label{eq1_7}
0\longrightarrow HC_1(A)\longrightarrow \left(\Lambda^2
A\right)/{\tt ({\it relations}\ (\ref{eqrev9}))}\longrightarrow
[A,A]\longrightarrow 0
\end{equation}
where the last map is the commutator map: $a\wedge b\mapsto [a,b]$.
The correctness of this map follows from the following relation for
any associative algebra $A$:
$$
[a,bc]+[b,ca]+[c,ab]=0
$$
The kernel of this map is the first cyclic homology $HC^1(A)$. Now
the statement (i) of lemma follows from the fact that $HC^1(A_n)=0$.

Consider $A$ as a Lie algebra with the bracket $[a,b]=a\cdot
b-b\cdot a$. Then the short exact sequence (\ref{eq1_7}) is a
sequence of $A$-modules, and the action of $A$ on $HC^1(A)$ is
trivial. Consider the corresponding long exact sequence in Lie
algebra homology:
$$
\dots\rightarrow HC_1(A)\rightarrow
\left(\Lambda^2(A/[A,A])\right)/(a\wedge (b\cdot c)+b\wedge (c\cdot
a)+c\wedge(a\cdot b)=0)\rightarrow [A,A]/[A,[A,A]]\rightarrow 0
$$
One should explain the middle term. It is not true that
$\Lambda^2(A)/[A,\Lambda^2(A)]$ is $\Lambda^2(A/[A,A])$, and we
should use the identity $a\wedge (b\cdot c)+b\wedge (c\cdot
a)+c\wedge(a\cdot b)=0$. We have by this identity:
$$
[a,b\wedge c]=[a,b]\wedge c+b\wedge [a,c]=a\wedge [b,c]
$$
which explains the middle term. Again, the statement (ii) of lemma
follows from the fact that $HC^1(A_n)=0$.
\end{proof}
\end{lemma}
\begin{remark}
The definition of a $k$-form on an associative algebra $A$,
generalizing the definition of 1-form here, is given in [K2]. It is
proven there that the de Rham complex of the algebra $A_n$ obeys the
Poincare lemma. It is also proven that the space of closed two-forms
on $A_n$ is $[A_n,A_n]$. From this point of view, the statement (i)
of the Lemma above is the Poincare lemma at the first term.
\end{remark}
\subsection{The case $n=2$}
We denote $A_{n,2}=[A_n,A_n]/[A_n,[A_n,A_n]]$, and we denote by
$A_{n,2}^{\ell}$ the graded component in $A_{n,2}$ consisting from
the monomials of the length $\ell$. We prove here the following
theorem:
\begin{theorem*}
$\dim A_{2,2}^\ell=\ell-1$
\begin{proof}
By Lemma 1.2 we know that $[A_n,A_n]/[A_n,[A_n,A_n]]\simeq
\left(\Lambda^2(A_n/[A_n,A_n])\right)/({\it relations}\
(\ref{eqrev9}))$. The idea is to show that any element in
$\left(\Lambda^2(A_2/[A_2,A_2])\right)/({\it relations}\
(\ref{eqrev9}))$ is equivalent to an element of the form
$x_1^k\wedge x_2^m$, where $k,m\ge 1$.

On the other hand, the space $\Omega^2_{closed}(\mathbb{C}^2)$ of
closed 2-forms of "length" $\ell$ (the length is the eigenvalue of
the operator $Lie_e$ where $e=x_1\frac{\partial}{\partial
x_1}+x_2\frac{\partial}{\partial x_2}$ is the Euler vector field) on
the 2-dimensional vector space has also dimension $\ell-1$. Denote
by $\wtilde{m}$ the image of a monomial $m\in A_2$ to the
commutative algebra $\mathbb{C}[x_1,x_2]=A_2/A_2[A_2,A_2]A_2$. Then
we have the map $[m_1,m_2]\mapsto d(\wtilde{m_1})\wedge
d(\wtilde{m_2})$, which clearly defines a map
$[A_2,A_2]/[A_2,[A_2,A_2]]\to\Omega^2_{closed}(\mathbb{C}^2)$. We
want to prove that this map is an isomorphism.

We start with the isomorphism $\Omega^1_{A_2}/d\Omega^0_{A_2}\simeq
[A_2,A_2]$. Now we consider the space of non-commutative 1-forms
$adb$ modulo exact forms and modulo forms of the types $a\cdot
d([x,y])$ and $[a,b]dx$. This space is isomorphic to the quotient
$[A_2,A_2]/[A_2,[A_2,A_2]]$. Let us compute this space.

Any 1-form can be reduced to the sum of forms of the type
$(dx_1)a(x_1,x_2)$ and $d(x_2)b(x_1,x_2)$ where $a$ and $b$ are
non-commutative monomials. Consider the form $(dx_1)a(x_1,x_2)$.
Suppose $a=a_1\cdot a_2$, then in our quotient space
\begin{equation}\label{eqrev10}
(dx_1)\cdot a_1\cdot a_2=(dx_1)\cdot a_2\cdot a_1
\end{equation}
Using this operation, we can suppose that $a(x_1,x_2)$ has form
$a=x_1^k\cdot a_1$ where $a_1$ starts with $x_2$ (or is equal to 1).
We have in our quotient space:
\begin{equation}\label{eqrev11}
(dx_1)\cdot x_1^k\cdot a_1=\frac1{k+1} (dx_1^{k+1})\cdot a_1
\end{equation}
Indeed, $(dx_1^{k+1})\cdot a_1=\sum_{i=0}^k x_1^i \cdot dx_1\cdot
x_1^{k-i}\cdot a_1$. Consider a summand $x_1^i \cdot dx_1\cdot
x_1^{k-i}\cdot a_1$. Using the cyclic symmetry of 1-forms, the
latter is the same that $dx_1\cdot x_1^{k-i}\cdot a_1\cdot x_1^i$.
Now using the equation (\ref{eqrev10}), we see that it is the same
that $dx_1\cdot x_1^k\cdot a_1$. Equation (\ref{eqrev11}) is proven.

Now we proceed in the same way: we represent $a_1$ as
$a_1=x_2^m\cdot x_1^n\cdot a_2$. Using the property (\ref{eqrev10})
we see that $dx_1\cdot x_2^m\cdot x_1^n\cdot a_2=dx_1\cdot
x_1^n\cdot a_2\cdot x_2^m=\frac1n d(x_1^{k+1+n})\cdot (a_2\cdot
x_2^m)$, and so on. Finally, we obtain that $(dx_1)\cdot a(x_1,
x_2)$ is equivalent to a form of the type $d(x_1^N)\cdot x_2^M$,
which proves the theorem.

\comment We prove it in two steps using the relation $a\wedge
(b\cdot c)+b\wedge (c\cdot a)+c\wedge(a\cdot b)=0$.

{\it Step 1}. We prove that a general element in
$\left(\wedge^2(A_n/[A_n,A_n])\right)/({\tt relations\ (1)})$ of the
form $M_1\wedge M_2$ is equivalent to an element of the form
$x^k\wedge M_3$ (here $M_1, M_2, M_3$ are noncommutative {\it
cyclic} words in $x$ and $y$).

Represent the cyclic word $M_1$ as $M_1=A\cdot B$ where
$A=(yyy...yy)$ and $B$ starts and ends in $x$'s, and represent
$M_2$ as $C=(yyy...y)\cdot\wtilde{M_2}\cdot (x...x)$. It is
possible to do if the both words $M_1$ and $M_2$ contain $y$'s,
otherwise we are done. Define the number of cuts $c_{M_i}$ in each
word $M_i$ as the number of places in the {\it cyclic} word where
$x$ bounds with $y$. Define the total number of cuts $c_{M_1,
M_2}=c_{M_1}+c_{M_2}$. Now we use the identity $(1)$:
\begin{equation}\label{eq1_12}
M_1\wedge M_2=-(B\cdot C)\wedge A-(C\cdot A)\wedge B
\end{equation}
Now it is clear that the total number of cuts $c_{B\cdot C,A}$ and
$c_{C\cdot A,B}$ are {\it less} than $c_{M_1,M_2}$. Then we can
iterate this process, and finally we obtain an equivalent monomial
of the form $x^k\wedge M_3$.

{\it Step 2}. We prove that each monomial of the form $x^k\wedge
M_3$ is equivalent to a monomial of the form $x^k\wedge y^m$.

If we have a monomial of the form $x^k\wedge M_3$ we can not
proceed as above. Our goal is to move all $x$'s from $M_3$ to the
left factor $x^k$. We find a maximal string of $x$'s in $M_3$:
$M_3=x^s\cdot M_4$ where $M_4$ starts and ends in $y$'s. Next, we
can suppose that $s\ge k$. Otherwise we can write by the identity
$a\wedge (b\cdot c)+b\wedge (c\cdot a)+c\wedge(a\cdot b)=0$:
\begin{equation}\label{eq1_13}
x^k\wedge (x^s\cdot M_4)=x^s\wedge (M_4\cdot x^k)-x^{k+s}\wedge
M_4
\end{equation}
The last summand has strictly less number of cuts, and the first
one is equal to $x^s\wedge (x^k\cdot M_4)$ because $M_4\cdot x^k$
here is considered as a cyclic word. Therefore we can suppose that
$s\ge k$.

Now put $A=x^k, B=x^k, C=x^{s-k}\cdot M_4$. Then we get from our
equation:
$$
x^k\wedge (x^k\cdot (x^{s-k}\cdot M_4))=x^{2k}\wedge (x^{s-k}\cdot
M_4)-x^k\wedge (x^k\cdot (x^{s-k}\cdot M_4))
$$
that is
\begin{equation}\label{eq1_14}
\frac12 x^k\wedge (x^s\cdot M_4)=x^{2k}\wedge (x^{s-k}\cdot M_4)
\end{equation}
It is clear that using the equations (\ref{eq1_13}) and
(\ref{eq1_14}) we can consequently move all $x$'s to the left
factor.

We have proved that each monomial $M_1\wedge M_2$ is equivalent to
$x^k\wedge y^m$.

{\it Step 3}. We have proved that $\dim [A_2,A_2]/[A_2,[A_2,A_2]]\le
\ell-1$. Now we need to prove that this dimension is exactly
$\ell-1$. It is actually clear from our description of
$[A_2,A_2]/[A_2,[A_2,A_2]]$ as
$\left(\wedge^2(A_n/[A_n,A_n])\right)/({\tt relations\
(\ref{eqrev9})})$. Indeed, we can not obtain any linear equations on
$x^k\wedge y^s$, $k+s=\ell$, using the relations $a\wedge (b\cdot
c)+b\wedge (c\cdot a)+c\wedge(a\cdot b)=0$.
\endcomment
\end{proof}
\end{theorem*}

\subsection{The case $n=3$}
\begin{theorem*}
The dimension of the space $A_{3,2}^\ell$ is $\ell^2-1$.
\begin{proof}
Again, it is the dimension of the polynomial closed two-forms on
$\mathbb{C}^3$ of the length $\ell$. Our proof is analogous to the
proof in the case $n=2$.

First we reduce a non-commutative 1-form to the form $d(x_1^k)\cdot
a(x_2, x_3)$. Then, modulo exact forms, it is $-da(x_2,x_3)\cdot
x_1^k$. We can proceed as above to reduce any form to the type
$d(x_i^{k_1})\cdot x_j^{k_2}\cdot x_m^{k_3}$ where $x_j^{k_2}$ and
$x_m^{k_3}$ commute ($\{i,j,m\}=\{1,2,3\}$).
\end{proof}
\end{theorem*}
We can not apply this proof for $n>3$. On the other hand,
computations showed that, starting from $n=4$, the dimension of the
space $A_{n,2}^\ell$ is {\it greater} than the dimension of the
corresponding closed 2-forms on $\mathbb{C}^n$ (see (17) in Section
4.2). We give the answer in the next Section.

\section{The quotient $[A_n,A_n]/[A_n,[A_n,A_n]]$ for general $n$}
Consider the following algebra structure on the (commutative) even
forms on an $n$-dimensional vector space
$\Omega^{even}(\mathbb{C}^n)$:
\begin{equation}\label{eqrev1}
\omega_1\circ\omega_2=\omega_1\wedge\omega_2+d\omega_1\wedge
d\omega_2
\end{equation}
where $d$ is the de Rham differential. Notice that this product is
not commutative neither skew-commutative. Later on, we consider only
this algebra structure on $\Omega^{even}$.
\begin{remark}
For any (not necessarily commutative) differential graded
associative algebra $A^\mb$ we can define a new algebra
$A^\mb_\star$ with the product \begin{equation}\label{eqnew1} a\star
b=a\cdot b+(-1)^{\deg a}(da)\cdot (db)
\end{equation}
($a,b$ are homogeneous) which is also associative.
\end{remark}

We have a map $\varphi_n\colon A_n\to\Omega^{even}(\mathbb{C}^n)_*$
which maps $x_i\in A_n$ to $x_i\in \Omega^0(\mathbb{C}^n)$, and we
extend it to $A_n$ in the unique way to get a map of algebras.
\subsection{}
\subsubsection{}
\begin{lemma}
The map $\varphi_n\colon A_n\to \Omega^{even}(\mathbb{C}^n)_*$ is
surjective.
\begin{proof}
Prove first that any monomial on the coordinates $\{x_i\}$'s belongs
to the image. Let $M_1$ and $M_2$ be two such monomials which belong
to the image of $\varphi_n$, we prove that $M_1\cdot M_2$ also does.
If $M_1=\varphi_n(R_1)$, and $M_2=\varphi_n(R_2)$, then $M_1\cdot
M_2=\frac12\varphi_n(R_1\cdot R_2+R_2\cdot R_1)$. It proves that any
monomial on $\{x_i\}$'s belongs to the image because linear
monomials $x_1,\dots,x_n$ belong to the image by definition.
Analogously we prove that any even monomial on $\{dx_i\}$'s belongs
to the image, and the general statement.
\end{proof}
\end{lemma}
\subsubsection{}
\begin{lemma}
\begin{itemize}
\item[(i)] The map $\varphi_n$ maps the commutator
$[A_n,A_n]$ to the closed (=exact) forms of degree $>0$
$\Omega^{even+}_{closed}(\mathbb{C}^n)$, and the map
$\varphi_n\colon [A_n,A_n]\to\Omega^{even+}_{closed}(\mathbb{C}^n)$
is surjective,
\item[(ii)] the triple commutator
$[A_n,[A_n,A_n]]$ is mapped  by $\varphi_n$ to 0,
\item[(iii)] the kernel of the map $\varphi_n$ is
$K_n=A_n\cdot[A_n,[A_n,A_n]]\cdot A_n$.
\end{itemize}
\begin{proof}
(i): it is clear that $[\omega_1,\omega_2]=2d\omega_1\wedge
d\omega_2$. Therefore, the image $[A_n,A_n]$ belongs to closed
forms. Surjectivity can be proved analogously with the lemma above,

(ii) it is clear from (i),

(iii) it follows from (ii) that $K_n$ belongs to the kernel of
$\varphi_n$, because $\varphi_n$ is a map of associative algebras,
and its kernel is a two-sided ideal. Now it is sufficiently to prove
that the algebra $A_n/K_n$ is isomorphic under $\varphi_n$ to
$\Omega^{even}(\mathbb{C}^n)_*$. It follows from the following
presentation of the {\it commutative} algebra
$\Omega^{even}(\mathbb{C}^n)$ by generators and relations: it is
generated by $\{x_i\}$ and $\{dx_i\wedge dx_j\}$ with the usual
commutativity relations and the relation
$$
(dx_i\wedge dx_j)\cdot (dx_k\wedge dx_l)=-(dx_i\wedge dx_k)\cdot
(dx_j\wedge dx_l)
$$
Therefore, $\Omega^{even}(\mathbb{C}^n)$ is a commutative algebra
generated by $\{x_i\}$ and $\{\eta_{i,j}\}$, where
$\eta_{i,j}=-\eta_{j,i}$ and with the relations
\begin{equation}\label{strange1}
\eta_{i,j}\cdot \eta_{k,l}+\eta_{i,k}\cdot\eta_{j,l}=0
\end{equation}
Now let us consider the algebra $A_n/K_n$. It is generated by
$\{x_i\}$ and $\{[x_i,x_j]\}$. We should check the relations
(\ref{strange1}), that is
\begin{equation}\label{strange2}
[x_i,x_j]\cdot [x_k,x_l]+[x_i,x_k]\cdot[x_j,x_l] \in K_n
\end{equation}
It follows from the following identity in the free algebra:
\begin{equation}\label{strange3}
\begin{aligned}
\ & [x_i,x_j]\cdot [x_k,x_l]+[x_i,x_k]\cdot [x_j,x_l]=\\
&[[x_j,x_k],x_ix_l]+x_i[x_k,[x_j,x_l]]+[[x_i,x_j],x_k]x_l-[[x_ix_l,x_k],x_j]
\end{aligned}
\end{equation}
Lemma is proven.
\end{proof}
\end{lemma}
\subsection{The main theorem}
We prove here the following theorem:
\begin{theorem*}
The map $\varphi_n$ induces an isomorphism $\varphi_n\colon
[A_n,A_n]/[A_n,[A_n,A_n]]\eqto\Omega^{even+}_{closed}$.
\end{theorem*}
By Lemma 2.1.2 above, the Theorem follows from the following Lemma:
\begin{klemma}
The intersection $[A_n,A_n]\bigcap\left( A_n\cdot[A_n,[A_n,A_n]]
\cdot A_n\right)=[A_n,[A_n,A_n]]$.
\end{klemma}
We prove this Lemma and the Theorem in the rest of this Section and
in Section 3.
\subsubsection{}
Consider the space $[A_n,A_n\cdot [A_n,[A_n,A_n]]]$. It is clear
that this space belongs to the intersection $[A_n,A_n]\bigcap\left(
A_n\cdot[A_n,[A_n,A_n]] \cdot A_n\right)$. We first prove the
Key-Lemma in this particular case.
\begin{lemma}
The space $[A_n,A_n\cdot [A_n,[A_n,A_n]]]$ belongs to
$[A_n,[A_n,A_n]]$.
\begin{proof}
Let $t_1,t_2,t_3,t_4,t_5$ be arbitrary monomials in $A_n$. We need
to prove that
\begin{equation}\label{eqrev2}
[t_1,t_2\cdot [t_3,[t_4,t_5]]]\in [A_n,[A_n,A_n]]
\end{equation}
We have:
\begin{equation}\label{eqrev3}
\begin{aligned}
\ [t_1,t_2[t_3,[t_4,t_5]]]&=[t_1,t_2t_3[t_4,t_5]]-[t_1,
t_2[t_4,t_5]t_3]\\
&=[t_1,t_2t_3[t_4,t_5]]-[t_1,
t_3t_2[t_4,t_5]]+[t_1,[t_3,t_2[t_4,t_5]]]
\end{aligned}
\end{equation}
Notice that the third summand in the last line belongs to
$[A_n,[A_n,A_n]]$. Now we apply the identity:
\begin{equation}\label{eqrev4}
[a,bc]+[b,ca]+[c,ab]=0
\end{equation}
We set $a=t_1$, $b=t_2t_3$ or $t_3t_2$, and $c=[t_4,t_5]$. By
(\ref{eqrev4}) we have:
\begin{equation}\label{eqrev5}
\begin{aligned}
\ [t_1,t_2[t_3,[t_4,t_5]]]&=[t_1,t_2t_3[t_4,t_5]]-[t_1,
t_3t_2[t_4,t_5]]+[t_1,[t_3,t_2[t_4,t_5]]]\\
&=-[t_2t_3,[t_4,t_5]t_1]-[[t_4,t_5], t_1t_2t_3]\\
&+[t_3t_2,[t_4,t_5]t_1]+[[t_4,t_5], t_1t_3t_2]\\
&+[t_1,[t_3,t_2[t_4,t_5]]]\\
&=[[t_3,t_2],[t_4,t_5]t_1]-[[t_4,t_5], t_1t_2t_3]+[[t_4,t_5],
t_1t_3t_2]+[t_1,[t_3,t_2[t_4,t_5]]]
\end{aligned}
\end{equation}
\end{proof}
\end{lemma}
\subsubsection{The proof of the Theorem}
By Lemma 1.2, we have the isomorphism
$\theta\colon\left(\Lambda^2(A_n/[A_n,A_n])\right)/({\it relations}\
(2))\eqto[A_n,A_n]/[A_n,[A_n,A_n]]$. The map $\theta$ is induced by
the map $\theta\colon a\wedge b\mapsto [a,b]$.

Recall that we denote by $K_n$ the kernel of the map of algebras
$\varphi_n\colon A_n\to \Omega^{even}(\mathbb{C}^n)_*$, that is, the
space $K_n=A_n\cdot[A_n,[A_n,A_n]]\cdot A_n$ by Lemma 2.1.2(iii). By
Lemma 2.2.1, the bracket $[K_n,A_n]\in[A_n,[A_n,A_n]]$, and,
therefore, the map $\theta$ defines a map
\begin{equation}\label{eqrev7}
\theta\colon\left(\Lambda^2(A_n/(K_n+[A_n,A_n]))\right)/({\it
relations}\ (2))\eqto[A_n,A_n]/[A_n,[A_n,A_n]]
\end{equation}
Here in the last formula we reduce the relation $a\wedge bc+b\wedge
ca+c\wedge ab=0$ modulo $K_n+[A_n,A_n]$. \subsubsubsection{}
\begin{lemma}
The space $A_n/(K_n+[A_n,A_n])$ is isomorphic to
$\Omega^{even}(\mathbb{C}^n)/\mathrm{Im}d$, and the isomorphism is
given by the map $\varphi_n$.
\begin{proof}
It follows from Lemma 2.1.1 and 2.1.2
\end{proof}
\end{lemma}
\subsubsubsection{} Now we deduce our Theorem to the following
result:
\begin{lemma}
The space $\Lambda^2\left(\Omega^{even}/\mathrm{Im}d\right)/({\it
relations}\ \omega_1\bigwedge(\omega_2\wedge\omega_3)+
\omega_2\bigwedge(\omega_3\wedge\omega_1)+\omega_3\bigwedge(\omega_1\wedge\omega_2)=0)
$ is isomorphic to $\Omega^{even+}_{closed}$, and the isomorphism is
given by the formula $\alpha\bigwedge\beta\mapsto d\alpha\wedge
d\beta$.
\end{lemma}

We prove Lemma 2.2.2.2 in the next Section.

\section{A proof of Lemma 2.2.2.2}
\subsection{The irreducible $W_n$-modules of the class
$\mathcal{C}$} We need to prove that that the space
$\Lambda^2\left(\Omega^{even}/\mathrm{Im}d\right)/({\it relations}\
\omega_1\bigwedge(\omega_2\wedge\omega_3)+
\omega_2\bigwedge(\omega_3\wedge\omega_1)+\omega_3\bigwedge(\omega_1\wedge\omega_2)=0)
$ is isomorphic to $\Omega^{even+}_{closed}$. The both sides are
modules over the Lie algebra $W_n$ of polynomial vector fields on an
$n$-dimensional vector space $\mathbb{C}^n$. We are going to define
some (a very general) class $\mathcal{C}$ of $W_n$-modules; all
$W_n$-modules we meet here belong to this class.

Consider the vector field $e=\sum_{i=1}^n
x_i\frac{\partial}{\partial x_i}$. We characterize the class
$\mathcal{C}$ of $W_n$-modules $L$ by the two conditions:
\begin{itemize}
\item[(i)] the operator $e$ is semisimple on $L$ with
finite-dimensional eigenspaces,
\item[(ii)] the eigenvalues of $e$ are bounded from below on $L$.
\end{itemize}

We are going to describe all irreducible modules over $W_n$ of the
class $\mathcal{C}$. Let $W_n^0$ be the Lie subalgebra of $W_n$ of
vector fields vanishing at the origin. Then $W_n^0$ has the
subalgebra $W_n^{00}$ of vector fields vanishing at the origin with
zero of at least second order. Actually $W_n^{00}$ is an ideal in
$W_n^0$, and $W_n^0/W_n^{00}\eqto\gl_n$.

Let $D$ be a Young diagram, and let $F_D$ be the corresponding
$\gl_n$-module (see [Ful]). Denote by $\mathcal{F}_D$ the coinduced
module $\mathcal{F}_D=Hom_{U(W_n^0)}(U(W_n),F_D)$.

\begin{theorem*}
\begin{itemize}
\item[(i)] The all representations $\mathcal{F}_D$ are irreducible
except the case when $D$ is just a one column, that is
$F_D=\Lambda^i(\mathbb{C}^n)^*$,
$\mathcal{F}_D=\Omega^i(\mathbb{C}^n)$; in the last case
$\mathcal{F}_D$ contains the image of the de Rham differential
$d\Omega^{i-1}(\mathbb{C}^n)$, which is irreducible,
\item[(ii)] the modules $\mathcal{F}_D$ for $D$ not a 1 column,
$\Omega^i(\mathbb{C}^n)/d\Omega^{i-1}(\mathbb{C}^n)$, and the
trivial representation exhaust all irreducible $W_n$-modules of the
class $\mathcal{C}$.
\end{itemize}
\end{theorem*}
\subsection{}
Let $V$ be an $n$-dimensional vector space, we prefer to work in
not-coordinate way. Then $\Omega^i(V)$ is coinduced from
$\Lambda^i(V^*)$.
\begin{lemma}
The map $i\colon V^*\otimes
\left(\Omega^{even}(V)/\mathrm{Im}d\right)\to
\Lambda^2\left(\Omega^{even}/\mathrm{Im}d\right)/({\it relations}\
\omega_1\bigwedge(\omega_2\wedge\omega_3)+
\omega_2\bigwedge(\omega_3\wedge\omega_1)+\omega_3\bigwedge(\omega_1\wedge\omega_2)=0)$
is surjective. Here we consider $V^*$ as linear 0-forms, and the map
$i$ is the composition of the inclusion with the subsequent
factorization.
\begin{proof}
Denote by $\overline{\omega_1\wedge\omega_2}$ the class of
$\omega_1\wedge\omega_2$ in the quotient space
$\Lambda^2\left(\Omega^{even}/\mathrm{Im}d\right)/({\it relations}\
\omega_1\bigwedge(\omega_2\wedge\omega_3)+
\omega_2\bigwedge(\omega_3\wedge\omega_1)+\omega_3\bigwedge(\omega_1\wedge\omega_2)=0)$.
Write $\omega_1=x_i\wedge\omega_1^{(1)}$. We have:
$x_i\omega_1^{(1)}\bigwedge\omega_2+(\omega_1^{(1)}\wedge\omega_2)\bigwedge
x_i+(\omega_2 x_i)\bigwedge\omega_1^{(1)}=0$ in the quotient space.
The second summand $(\omega_1^{(1)}\wedge\omega_2)\bigwedge x_i$
belongs to $\overline{V^*\otimes(\Omega^{even}(V)/\mathrm{Im}d)}$.
Thus, modulo this image, we can freely move all $x_i$'s from
$\omega_1$ to $\omega_2$. We can do it many times. Finally, instead
of $\omega_1$ we will have a form without $x_i$'s, that is, a form
$dx_{i_1}\wedge\dots\wedge dx_{i_k}$. This form is exact, and
therefore is zero in
$\Lambda^2\left(\Omega^{even}/\mathrm{Im}d\right)$. We are done.
\end{proof}
\end{lemma}
\subsection{}
Thus, we have a surjective map $i$ of $V\otimes
\Omega^{2k}(V)/\mathrm{Im}d$ to
$\Lambda^2\left(\Omega^{even}/\mathrm{Im}d\right)/({\it relations}\
\omega_1\bigwedge(\omega_2\wedge\omega_3)+
\omega_2\bigwedge(\omega_3\wedge\omega_1)+\omega_3\bigwedge(\omega_1\wedge\omega_2)=0)$,
and the latter is mapped to $\Omega^{2k+2}_{closed}(V)$ by the
formula $\omega_1\wedge\omega_2\mapsto d\omega_1\wedge d\omega_2$.
This map $j$ is clearly surjective. We need to prove that the
induced map $(V\otimes \Omega^{2k}(V)/\mathrm{Im}d)/{\it
relations}\to\Omega^{2k+2}_{closed}(V)$ is an isomorphism. Our tool
is Theorem 4.1.

The main point is that $(V\otimes \Omega^{2k}(V)/\mathrm{Im}d)$ is
{\it not} a $W_n$-module, only  $(V\otimes
\Omega^{2k}(V)/\mathrm{Im}d)/{\it relations}$ is. But it is still a
$\gl_n$-module. Therefore, we should use the representation theory
of $\gl_n$-modules.

First of all, we describe the $\Omega^k(V)$ as a $\gl_n$-module. The
answer is given in Figure 1. It follows from the
Littlewood-Richardson rule applied to $\Lambda^k(V^*)\otimes
S^N(V^*)$. \sevafigc{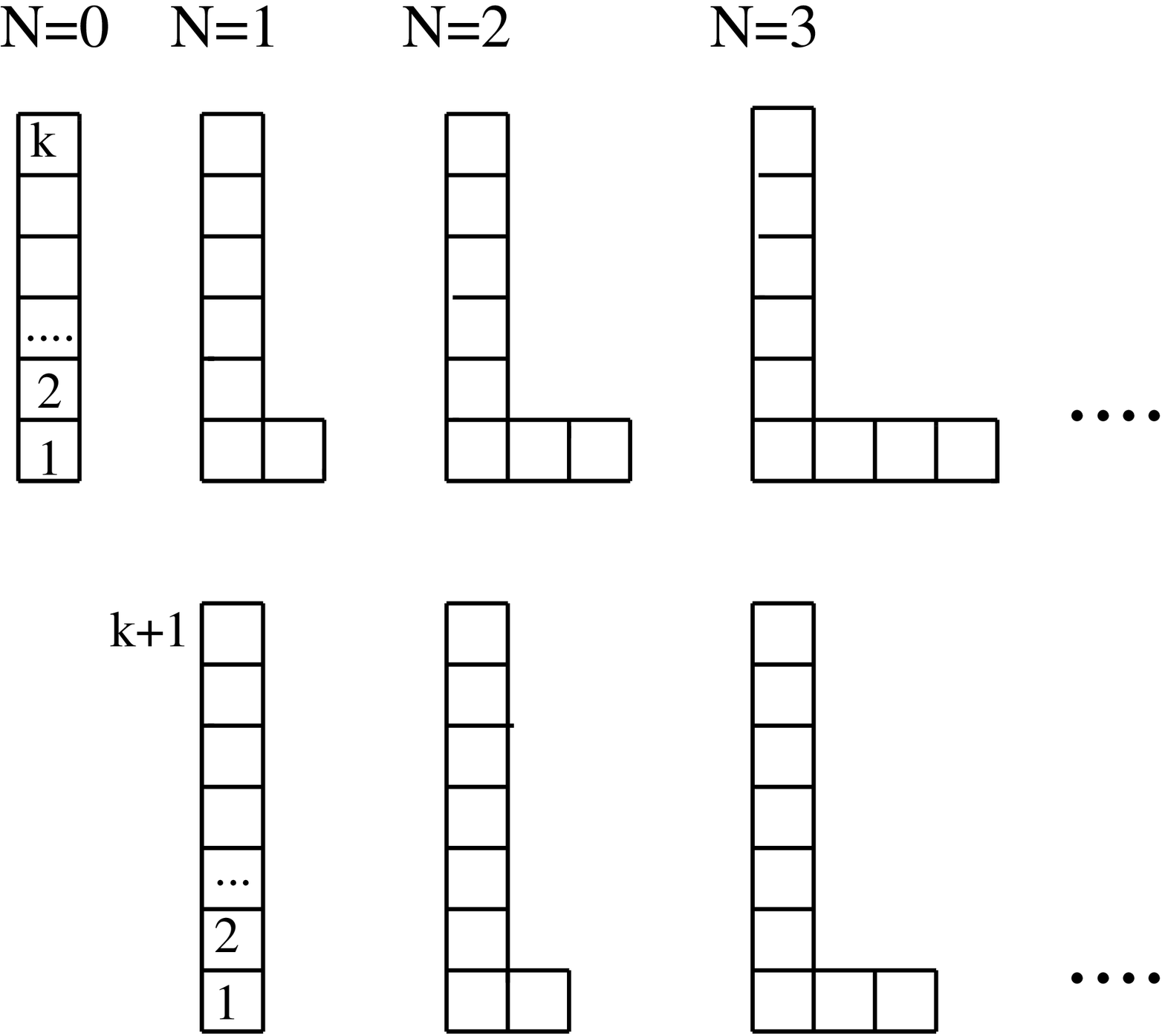}{100mm}{0}{The space
$\Omega^k(V)$ as a $\gl_n$-module} As a $W_n$-module, $\Omega^k(V)$
has a submodule $Im_k=\Omega^k(V)_{closed}=d\Omega^{k-1}(V)$, and
the quotient module $\Omega^k(V)/Im_k\simeq Im_{k+1}$ by the
Poincare lemma. It is clear that as a $\gl_n$-module, the submodule
$Im_k$ is the first line in Figure 1, while the quotient-module
$\Omega^k(V)/Im_k\simeq Im_{k+1}$ is the second line.

Now we should pass from $\Omega^k(V)/d\Omega^{k-1}(V)$ to $V\otimes
(\Omega^k(V)/d\Omega^{k-1}(V))$. The space
$\Omega^k(V)/d\Omega^{k-1}(V)$ is the second line in Figure 1, and
now we add a 1 new box to each Young diagram by the
Littlewood-Richardson rule. We obtain a number of diagrams, but the
only one among them will be a column (of the height $k+2$).

Now consider the situation of Lemma 2.2.2.2 ($k$ is even, $k=2m$).
The only column of the height $2m+2$ maps isomorphically to the
corresponding component in $\Omega^{2m+2}(V)_{closed}$ under the map
of Lemma 2.2.2.2. Therefore, there are no columns in the kernel of
the map $i\colon V^*\otimes
\left(\Omega^{even}(V)/\mathrm{Im}d\right)\to
\Lambda^2\left(\Omega^{even}/\mathrm{Im}d\right)/({\it relations}\
\omega_1\bigwedge(\omega_2\wedge\omega_3)+
\omega_2\bigwedge(\omega_3\wedge\omega_1)+\omega_3\bigwedge(\omega_1\wedge\omega_2)=0)$.
But when we consider the map $\bar{i}\colon  (V^*\otimes
\left(\Omega^{even}(V)/\mathrm{Im}d\right))/{\it relations} \to
\Lambda^2\left(\Omega^{even}/\mathrm{Im}d\right)/({\it relations}\
\omega_1\bigwedge(\omega_2\wedge\omega_3)+
\omega_2\bigwedge(\omega_3\wedge\omega_1)+\omega_3\bigwedge(\omega_1\wedge\omega_2)=0)$,
then its kernel is smaller then the kernel of the map $i$, and is a
$W_n$-module. It does not contain any column in the
$\gl_n$-decomposition. We know from Theorem 3.1 the description of
all $W_n$-modules of the class $\mathcal{C}$. It could not be
$d\Omega^l(V)$ because it does not contain any column. Therefore, it
it the coinduced module with some Young diagram $D$ which is not a
column. It is clear by the Littlewood-Richardson rule that $D\otimes
S^\mb(V^*)$ is {\it bigger} than we have from the second line of
Figure 1 multiplied by $V^*$.

Lemma 2.2.2.2 is proven \qed

Theorem 2.2 and the Key-Lemma 2.2 are proven.\qed

\comment
\section{The second proof of Lemma 3.2.2.2}
\subsection{The idea}
Recall that we denote by $\Omega^{even}_*$ the associative algebra
of polynomial even degree differential forms on a vector space $V$
with the product
$\omega_1\star\omega_2=\omega_1\cdot\omega_2+(d\omega_1)\cdot(d\omega_2)$.
The space of all differential forms is the space of functions on
$W=V\oplus V[1]$. We can consider $d\wedge d$, where $d$ is the de
Rham differential, as a bivector field on $W$. It is clear that this
bivector is Poisoon, that is, $[d\wedge d,d\wedge d]=0$. One can ask
ourselves the question how to construct a deformation quantization
of the algebra of forms on $V$ along this bivector, the answer is
\begin{equation}\label{eq10_1}
\omega_1\star\omega_2=\omega_1\cdot\omega_2+(-1)^{\deg
\omega_1}(d\omega_1)\cdot (d\omega_2)
\end{equation}
and our formula for the star-product on $\Omega^{even}_*$ is the
restriction of the formula (\ref{eq10_1}) product to the even forms.

Now we have the following very simple Lemma:
\begin{lemma}
The equality $HC_1(\Omega^{even}_*)=0$ implies Lemma 3.2.2.2.
\begin{proof}
By Lemma 1.2.1, we have the following long exact sequence for an
associative algebra
\begin{equation}\label{eq10_2}
\dots\rightarrow HC_1(A)\rightarrow
\left(\Lambda^2(A/[A,A])\right)/(a\wedge (b\cdot c)+b\wedge (c\cdot
a)+c\wedge(a\cdot b)=0)\rightarrow [A,A]/[A,[A,A]]\rightarrow 0
\end{equation}
Therefore, if $HC_1(A)=0$, we have an isomorphism
$\left(\Lambda^2(A/[A,A])\right)/(a\wedge (b\cdot c)+b\wedge (c\cdot
a)+c\wedge(a\cdot b)=0)\simeq [A,A]/[A,[A,A]]$.

Suppose that $A=\Omega^{even}_*$. Then
$[A,A]=\Omega^{even+}=d\Omega^{odd}$, and $[A,[A,A]]=0$, by Lemmas
in Section 3.1.
\end{proof}
\end{lemma}
Thus, it is enough to prove that $HC_1(\Omega^{even}_*)=0$. This is
not very easy. We first compute the Hochschild homology of (not
deformed) even forms on a vector space $V$, $HH_\mb(\Omega^{even})$.
The answer is quite non-trivial, and the only way to get it we know
is to use that $\Omega^{even}(V)$ is Koszul algebra. It follows from
the following theorem on Veronese powers (see [PP], Section 3).
\subsubsection{}
\begin{lemma}
Suppose that a $\mathbb{Z}_+$-graded quadratic Koszul algebra $A$ is
generated by $A_1$ and $A_0=\mathbb{C}$. Then for any $d\ge 2$ the
algebra $A(d)=\oplus_{n\ge 0}A_{nd}$ is also quadratic and Koszul.
\end{lemma}
We recall some definitions on Koszul algebras in the next
Subsection. In particular, the Hochschild (co)homology of any Koszul
algebra can be computed using some Koszul resolution. In general it
is not enough to complete the computation of the Hochschild
homology, but when $A$ is a {\it commutative} Koszul algebra we
develop in the next Subsection a mini-theory which gives the
explicit answer for $HH_1(A)$, the first Hochschild homology. Then
the cyclic homology is by the Connes-Tsygan bicomplex the quotient
of $HH_1(A)$ by the image of the "de Rham differential" $D\colon
HH_0(A)\to HH_1(A)$. We also compute in the next Subsection
explicitly this de Rham differential in the case of commutative
Koszul algebras. This will complete the computation of the first
cyclic homology $HC_1(\Omega^{even})$.
\subsection{\tt Koszul algebras, Koszul complex, and Hochschild
homology}
\subsection{ Cyclic homology}
\subsection{ First Hochschild and cyclic homology of a
commutative quadratic Koszul algebra} \subsubsection{} Let $A$ be a
commutative quadratic Koszul algebra with the space of generators
$V$. Consider the complex from Section 5.2, computing the Hochschild
homology of $A$. It is
\begin{equation}\label{eqmy1}
\dots\rightarrow [A^{!*}]^3\otimes A\rightarrow [A^{!*}]^2\otimes
A\rightarrow V\otimes A\rightarrow A\rightarrow 0
\end{equation}
Here $V$ in the term $V\otimes A$ is $V=[A^{!*}]^1$. The
differential $d\colon V\otimes A\to A$ is $d(v\otimes a)=v\cdot
a-a\cdot v$ is $0$ in the case of a commutative Koszul algebra. Let
us compute the image of the differential $d\colon [A^{!*}]^2\to
V\otimes A$. First of all, $[A^{!*}]^2\simeq I$ where $I\subset
V\otimes V$ is the quadratic part of the ideal of the quadratic
algebra $A$. Indeed, it is $\left((V^*\otimes V^*)/I^!\right)^*$
which is the subspace in $V\otimes V$, annihilating the ideal $I^!$,
which is $I$ by definition of $I^!$. Thus, we have the differential
\begin{equation}\label{eqmy2}
d\colon I\otimes A\to V\otimes A
\end{equation}
which for an element $i=\sum_k x_k\otimes y_k\in V\otimes V$ is
\begin{equation}\label{eqmy3}
d(i\otimes a)=\sum_k(x_k\otimes(y_k\cdot a)+y_k\otimes (x_k\cdot a))
\end{equation}
Finally, we gen the following explicit answer for the first
Hochschild homology $HH_1(A)$:
\begin{equation}\label{eqmy4}
HH_1(A)=\left( (V\otimes A)/\{\sum_k(x_k\otimes(y_k\cdot
a)+y_k\otimes (x_k\cdot a),\  \sum_kx_k\otimes y_k\in I\}\right)
\end{equation}

Now we are going to compute the first cyclic homology $HC_1(A)$. For
any associative algebra $A$,
\begin{equation}\label{eqmy5}
HC_1(A)=HH_1(A)/{\mathrm Im}\{D\colon HH_0(A)\to HH_1(A)\}
\end{equation}
where $D$ is the "de Rham differential" obtained from the
Connes-Tsygan bicomplex (see Section 5.3). We have the following
explicit formula for $D$:
\begin{equation}\label{eqmy6}
.......
\end{equation}
To prove (\ref{eqmy5}), one should just consider the spectral
sequence, associated with the Connes-Tsygan bicomplex.

Surprisingly, we have the following explicit description of $D$ in
the case of commutative Koszul algebras:
\begin{lemma}
Let $A$ be a commutative quadratic Koszul algebra with the space of
generators $V$. Consider the identity $\P=\sum_i(v_i\otimes v_i^*)$
in $V\otimes V^*$. Represent $A=T(V)/\Im$ where $\Im$ is the
quadratic ideal generated by $I$, where $T(V)$ is the free algebra
generated by $V$. Then the de Rham differential $D\colon HH_0(A)\to
HH_1(A)$ is equal to the action $D(a)=\sum_i
v_i\otimes\left(v_i^*(\tilde{a})\right)$. Here we lift $a\in
A=T(V)/\Im$ to any preimage of it in $T(V)$. Then $v_i^*$ acts on
$T(V)$ as a derivation satisfying the Leibniz rule and acting on the
generators $V$ just by the contraction. This definition is correct,
that is $D(a)$ does not depend on the choice of $\tilde{a}$.
\begin{proof}
First of all, let us prove that the definition of $D(a)$ is correct.
That mean that $\sum_i v_i\otimes\left(v_i^*(\tilde{a})\right)=0$ in
$HH_1(A)$ for any $\tilde{a}$ of the form $\tilde{a}\in\Im$, the
quadratic ideal in $T(V)$, generated by $I$. Consider for example
the expression $\sum_i v_i\otimes\left(v_i^*(i\cdot a_1\right)$
where $i=\sum_k x_k\otimes y_k\in I$. Then clearly
\begin{equation}\label{eqmy7}
\sum_i v_i\otimes\left(v_i^*(i\cdot a_1)\right)= \sum_k\left(
x_k\otimes (y_k\cdot a_1)+y_k\otimes (x_k\cdot a_1)\right)+\sum_i
v_i\otimes (i\cdot v_i^*(a_1))
\end{equation}
The last summand belongs to $\Im$ and therefore is 0 in $A$. Now the
correctness is clear because $HH_1(A)=\left( (V\otimes
A)/\{\sum_k(x_k\otimes(y_k\cdot a)+y_k\otimes (x_k\cdot a),\
\sum_kx_k\otimes y_k\in I\}\right)$.
..........................................
\end{proof}
\end{lemma}
\subsubsection{}
Moreover, we have the following analog of the
Hochschild-Kostant-Rosenberg theorem for $HH_1(A)$.
\begin{lemma}
Consider the map $A\otimes A\to HH_1(A)$, such that $a_1\otimes
a_2\mapsto a_1\cdot D(a_2)\in \left( (V\otimes
A)/\{\sum_k(x_k\otimes(y_k\cdot a)+y_k\otimes (x_k\cdot a),\
\sum_kx_k\otimes y_k\in I\}\right)$. Then it is an isomorphism on
the first Hochschild homology computed from the bar-resolution with
the first Hochschild homology computed via the Koszul resolution.
\qed
\end{lemma}
This theory completes the computation of $HC_1(A)$ for any
commutative quadratic Koszul algebra $A$. In the next Subsection we
are going to apply it to $A=\Omega^{even}(\mathbb{C}^n)$.
\begin{remark}
It would be very interesting to understand is it possible to
generalize the theory developed here for higher Hochschild homology
of a commutative quadratic Koszul algebra. In particular, how to
generalize the Hochschild-Kostant-Rosenberg map from the previous
lemma.
\end{remark}

\subsection{ A computation of $HC_1(\Omega^{even})$}
We compute here $HC_1(\Omega^{even})$. We have a short exact
sequence
\begin{equation}\label{eqmy10}
0\rightarrow HC_1(\Omega^{even})\rightarrow
\Lambda^2(\Omega^{even})/(a\wedge (b\cdot c)+b\wedge (c\cdot
a)+c\wedge(a\cdot b)=0)\rightarrow
[\Omega^{even},\Omega^{even}]\rightarrow 0
\end{equation}
which exists for any associative algebra. Here the most right term
is 0 because the algebra $\Omega^{even}$ is commutative. Therefore,
{\it a priori} we have:
\begin{equation}\label{eqmy11}
HC_1(\Omega^{even})\eqto\Lambda^2(\Omega^{even})/(a\wedge (b\cdot
c)+b\wedge (c\cdot a)+c\wedge(a\cdot b)=0)
\end{equation}
In our computation we replace $\Lambda^2(\Omega^{even})/(a\wedge
(b\cdot c)+b\wedge (c\cdot a)+c\wedge(a\cdot b)=0)$ by $(V\otimes
\Omega^{even})/{\tt some\ relations}$. Here $V$ will be
\begin{equation}\label{eqmy12}
V=\{x_i\}_{1\le i\le n}\oplus \{dx_i\wedge dx_j\}_{i<j}
\end{equation}
and our goal is to find the relations by which we should factorize
$V\otimes \Omega^{even}$ to get $\Lambda^2(\Omega^{even})/(a\wedge
(b\cdot c)+b\wedge (c\cdot a)+c\wedge(a\cdot b)=0)$.

By (\ref{eqmy5}) we have:
$$
HC_1(A)=HH_1(A)/{\mathrm Im}\{D\colon HH_0(A)\to HH_1(A)\}
$$
For a commutative Koszul algebra we know explicitly $HH_0(A)$,
$HH_1(A)$, and the de Rham operator $D$ from previous Subsection. We
know that the algebra $\Omega^{even}$ of even differential forms on
a vector space is Koszul with the generators $V$ (\ref{eqmy12}).
Denote $V=V_1\oplus V_2$ where $V_1=\{x_i\}$, and $V_2=\{dx_i\wedge
dx_j\}$.

We have: $HH_0(\Omega^{even})=\Omega^{even}$,
\begin{equation}\label{eqmy13}
\begin{aligned}
\ & HH_1(\Omega^{even})=[ (V\otimes
\Omega^{even})/(\sum_{i,j,k,l}(dx_i\wedge dx_j)\otimes(dx_k\wedge
dx_l\cdot a)+\\
& (dx_k\wedge dx_l)\otimes (dx_i\wedge dx_j\cdot a)+(dx_i\wedge
dx_k)\otimes(dx_j\wedge dx_l\cdot a)+(dx_j\wedge dx_l)\otimes
(dx_i\wedge dx_k\cdot a))]
\end{aligned}
\end{equation}
by (\ref{eqmy4}), because the only relations on $V$ are the
commutativity relations (which do not contribute to (\ref{eqmy4})),
and the relation
\begin{equation}\label{eqmy14}
(dx_i\wedge dx_j)\otimes (dx_k\wedge dx_l)+(dx_i\wedge dx_k)\otimes
(dx_j\wedge dx_l)=0
\end{equation}
(we suppose that $dx_i\wedge dx_j$ is skew-commutative in $i,j$).

Now we should factorize $HH_1(\Omega^{even})$ by the image of the de
Rham operator $D$. We use the explicit description of $D$ for a
commutative Koszul algebra from Lemma 5.4.1. We have:
\begin{equation}\label{eqmy15}
D(\omega)=\sum_ix_i\otimes \frac{\partial}{\partial x_i}\omega+
\sum_{i<j}(dx_i\wedge dx_j)\otimes i_{dx_i\wedge dx_j}\omega
\end{equation}
The first summand belongs to $V_1\otimes\Omega^{even}$, and the
second one belongs to $V_2\otimes\Omega^{even}$. Therefore, for a
homogeneous $\omega\in \Omega^{even}$, the equation $D\omega=0$ is
actually two equations:
\begin{equation}\label{eqmy16}
\sum_ix_i\otimes \frac{\partial}{\partial x_i}\omega=0
\end{equation}
and
\begin{equation}\label{eqmy17}
\sum_{i<j}(dx_i\wedge dx_j)\otimes i_{dx_i\wedge dx_j}\omega
\end{equation}
We had before
\begin{equation}\label{eqmy19}
\begin{aligned}
\ &\sum_{i,j,k,l}(dx_i\wedge dx_j)\otimes(dx_k\wedge dx_l\cdot
\omega)+(dx_k\wedge dx_l)\otimes (dx_i\wedge dx_j\cdot
\omega)+\\
&(dx_i\wedge dx_k)\otimes(dx_j\wedge dx_l\cdot \omega)+(dx_j\wedge
dx_l)\otimes (dx_i\wedge dx_k\cdot \omega)
\end{aligned}
\end{equation}
Finally, we get:
\begin{equation}\label{eqmy18}
HC_1(\Omega^{even})=V\otimes \Omega^{even}/({\tt relations\
(\ref{eqmy19}),\ (\ref{eqmy16}),\ (\ref{eqmy17})})
\end{equation}

It is not surprising from (\ref{eqmy11}) that we can not improve
this answer. Nevertheless, in the next Subsection we will factorize
$V\otimes \Omega^{even}$ by two extra equations, and after that this
{\it a priori} even more complicated expression will be nicely
become equal to $\Omega^{even+}_{closed}$.
\subsection{ A proof of Lemma 3.2.2.2}
Up to now we were computing $HC_1(\Omega^{even})$. It was explained
above in Section 5.1 that to prove Lemma 3.2.2.2 it is enough to
prove that $HC_1(\Omega^{even}_*)=0$ where $\Omega^{even}_*$ is the
deformed algebra (\ref{eq10_1}). We have the following short exact
sequence
\begin{equation}\label{eqmy20}
0\rightarrow HC_1(\Omega^{even}_*)\rightarrow
\Lambda^2(\Omega^{even})/(a\wedge (b* c)+b\wedge (c* a)+c\wedge(a*
b)=0)\rightarrow [\Omega_*^{even},\Omega_*^{even}]\rightarrow 0
\end{equation}
The most right term is
$[\Omega^{even}_*,\Omega^{even}_*]=\Omega^{even+}_{closed}$ by
Lemmas in Section 3.1. That is, to prove that
$HC_1(\Omega^{even}_*)=0$ is enough to prove that the last map is an
isomorphism
\begin{equation}\label{eqmy21}
\varphi\colon\Lambda^2(\Omega^{even})/(a\wedge (b* c)+b\wedge (c*
a)+c\wedge(a* b)=0)\eqto \Omega^{even+}_{closed}
\end{equation}
given by the map $\varphi\colon \omega_1\wedge\omega_2\mapsto
d\omega_1\wedge d\omega_2$.

The relation $a\wedge (b* c)+b\wedge (c* a)+c\wedge(a* b)=0$ is
equally two equations:
\begin{equation}\label{eqmy22}
\omega_1\bigwedge
(\omega_2\wedge\omega_3)+\omega_2\bigwedge(\omega_3\wedge
\omega_1)+\omega_3\bigwedge (\omega_1\wedge\omega_2)=0
\end{equation}
and
\begin{equation}\label{eqmy23}
\omega_1\bigwedge (d\omega_2\wedge
d\omega_3)+\omega_2\bigwedge(d\omega_3\wedge
d\omega_1)+\omega_3\bigwedge (d\omega_1\wedge d\omega_2)=0
\end{equation}
Notice now that $\Lambda^2(\Omega^{even})/(\ref{eqmy22})\eqto
HC_1(\Omega^{even})$ we have computed above using the Koszul
resolution. That is, to prove that $\varphi$ is an isomorphism it is
enough to prove that the corresponding to $\varphi$ map
\begin{equation}\label{eqmy24}
\varphi_1\colon HC_1(\Omega^{even})/(\ref{eqmy23})\eqto
\Omega^{even+}_{closed}
\end{equation}
is an isomorphism. We should do two things: to construct $\varphi_1$
which is $\omega_1\wedge\omega_2\mapsto d\omega_2\wedge d \omega_2$
but defined on $V\otimes \Omega^{even}/{\it relations})$, and to
factorize our answer for $HC_1(\Omega^{even})$ by the relation
(\ref{eqmy23}).

A problem here is that the relation (\ref{eqmy23}) is defined on
$\Lambda^2(\Omega^{even})$, not on $V\otimes \Omega^{even}$. For
this we have the "Hochschild-Kostant-Rosenberg  map" from Lemma
5.4.2. In its cyclic version, it is the map $\varphi_{HKR}\colon
\Lambda^2(\Omega^{even})\to HC_1(\Omega^{even})=(V\otimes
\Omega^{even})/({\tt relations\ (\ref{eqmy19}),\ (\ref{eqmy16}),\
(\ref{eqmy17})})$. This map is
$\varphi_{HKR}(\omega_1\bigwedge\omega_2)=\omega_1\cdot D(\omega_2)$
where $D$ is the de Rham differential. Therefore, we should
factorize $(V\otimes \Omega^{even})/({\tt relations\
(\ref{eqmy19}),\ (\ref{eqmy16}),\ (\ref{eqmy17})})$ by the relation
\begin{equation}\label{eqmy25}
\omega_1\cdot D(d\omega_2\wedge d\omega_3)+\omega_2\cdot
D(d\omega_3\wedge d\omega_1)+\omega_3\cdot D(d\omega_1\wedge
d\omega_2)
\end{equation}
The expression here belongs to $V\otimes\Omega^{even}$. Therefore it
gives two relations, because $V=V_1\oplus V_2$. These two relations
are:
\begin{equation}\label{eqmy26}
\sum_i\left(x_i\otimes(\omega_1\wedge \frac{\partial}{\partial
x_i}(d\omega_2\wedge d\omega_3))+x_i\otimes(\omega_2\wedge
\frac{\partial}{\partial x_i}(d\omega_3\wedge
d\omega_1))+x_i\otimes(\omega_3\wedge \frac{\partial}{\partial
x_i}(d\omega_1\wedge d\omega_2))\right)=0
\end{equation}
(for the component $V_1$), and
\begin{equation}\label{eqmy27}
\sum_{i<j}\left((dx_i\wedge dx_j)\otimes (\omega_1\wedge
[\frac{\partial}{\partial x_i}\omega_2\wedge
\frac{\partial}{\partial x_j}\omega_3-\frac{\partial}{\partial
x_j}\omega_2\wedge \frac{\partial}{\partial x_i}\omega_3]+{\tt
cycl}\right)
\end{equation}
(for the component $V_2$). Here {\tt cycl} means the sum of other
two terms obtaining from the permutations of indexes $(1,2,3)\mapsto
(2,3,1)$ and $(1,2,3)\mapsto (3,1,2)$.

Here we used the explicit formula for $D$ from (\ref{eqmy15}). Now
what is remained to prove is the following lemma:
\begin{lemma}
The map $\varphi_1\colon (V\otimes \Omega^{even})/({\tt relations\
(\ref{eqmy19}),\ (\ref{eqmy16}),\ (\ref{eqmy17}), (\ref{eqmy26}),
(\ref{eqmy27})})\to \Omega^{even+}_{closed}$ defined as 0 on
$V_2\otimes\Omega^{even}/{\it relations}$ and as
$x_i\otimes\omega\mapsto dx_i\wedge d\omega$ on $V_1\otimes
\Omega^{even}/{\it relations}$ is well defined, and is an
isomorphism.
\begin{proof}
\end{proof}
\end{lemma}

Lemma 3.2.2.2 is proven. \qed
\endcomment
\section{A $W_n$-action on $\tilde{gr}A_n$}
\subsection{The Theorem on $W_n$-action}
\comment We are going to formulate a Conjecture on all consecutive
quotients $A_{n,k}$.

The class $\mathcal{C}$ of $W_n$-modules introduced in Section 4 is
closed under extensions and subquotients.
\begin{conjecture}
All quotients $A_{n,k}$ for $k>1$ are $W_n$-modules of the class
$\mathcal{C}$. In particular, each quotient $A_{n,k}$ is a
$W_n$-module for $k>1$. The Lie bracket on
$gr^+(A_n)=\oplus_{k>1}A_{n,k}$ is $W_n$-equivariant.
\end{conjecture}

Moreover, the term $A_{n,1}=A_n/[A_n,A_n]$ contains the kernel $K_n$
of the map $\varphi_n\colon A_n\to\Omega^{even}$. The quotient
$A_{n,1}/K_n$ is isomorphic $\Omega^{even}/{\rm Im}d$ by Lemma
3.2.2.1. It is clear from the previous Section (Lemma 3.2.1) that
the kernel $K_n$ belongs to the {\it center} of the Lie algebra
$gr(A_n)$. Let us consider the Lie algebra
$\wtilde{gr}(A_n)=gr(A_n)/K_n$. Then this Lie algebra consists from
$W_n$-modules of the class $\mathcal{C}$, and the Lie bracket is
$W_n$-equivariant.

The first grading component of the Lie algebra
$\wtilde{gr}^1(A_n)=\Omega^{even}/{\rm Im}d$. The bracket
$\Lambda^2(\wtilde{gr}^1(A_n))\to \wtilde{gr}^2(A_n)=[A_n,A_n]/
[A_n,[A_n,A_n]]=\Omega^{even+}_{closed}$ is the map
$\omega_1\bigwedge\omega_2\to (d\omega_1)\wedge (d\omega_2)$ which
is clearly $W_n$-equivariant.

An evidence for the part of the Conjecture above stating that all
spaces $A_{n,k}$ are $W_n$-modules of the class $\mathcal{C}$ for
$k>1$ comes from the computation (\ref{eq1_3}). We have already
proved in the paper that $A_{n,2}$ is a $W_n$-module. Consider the
space $A_{n,3}$ for $n=2$. We are going to show that the dimensions
$\dim A_{2,3}^{\ell}$ are exactly like the character of a
$W_2$-module. Indeed, we know from (\ref{eq1_3}) that $\dim
A_{2,3}^\ell=2(\ell-2)$ for $3\le \ell\le 9$. Consider a
$W_2$-module coinduced from a Young diagram $D$ with a 2-dimensional
$\gl_2$-module on the level $\ell=3$. Then on a level $\ell$ this
$W_2$-module should have dimension $2 S_{\ell-3,2}$ where $S_{k,2}$
is dimension of symmetric polynomials in 2 variables of degree $k$.
We have: $S_{k,2}=k+1$. This shows that the spaces $A_{2,3}^\ell$
have the character of a $W_2$-module for $\ell\le 9$.

Consider now $A_{2,4}$. Here the polynomial for $A_{2,4}^\ell$ is
$3t^4+8t^5+13t^6+18t^7+23t^8+28t^9+\mathcal{O}(t^10)$. First take
the difference with $3t^4+6t^5+\dots+3(k+1)t^{k+4}+\dots$ which is
the character of a $W_2$-module. The difference is $\sum_{k\ge
0}2(k+1)t^{k+5}$ which is a character of $W_2$-module.

The case $A_{2,5}$ is analogous.

It is more interesting to consider the case $n=3$. Consider
$A_{3,3}$. \sevafigc{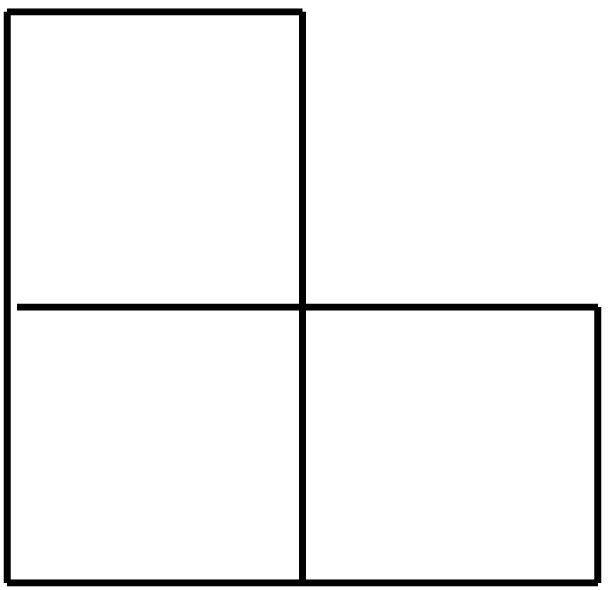}{40mm}{0}{The Young diagram $D_1$}
The polynomial for $A_{3,3}^\ell$ is
$p_{3,3}=8t^3+24t^4+48t^5+80t^6+\mathcal{O}(t^7)$. The first
coefficient $8$ is the space generated by Lie words
$[x_i,[x_j,x_k]]$ on 3 letters $x_1,x_2,x_3$ (some of $i,j,k$ may
coincide). This space has dimension 8, and as $\gl_3$-module it is
corresponded to the Young diagram $D_1$ (see Figure~2). The
irreducible $\gl_3$-module corresponded to this Young diagram has
dimension 8, and it can be easily computed by the character formula.
Consider the $W_3$-module coinduced from this $\gl_3$-module. This
coinduced character is $8\sum_{k\ge 0}\frac{(k+1)(k+2)}{2}t^{k+3}$.
It is exactly our $p_{3,3}$ up to $t^6$.

Finally, consider $A_{3,4}$. \sevafigc{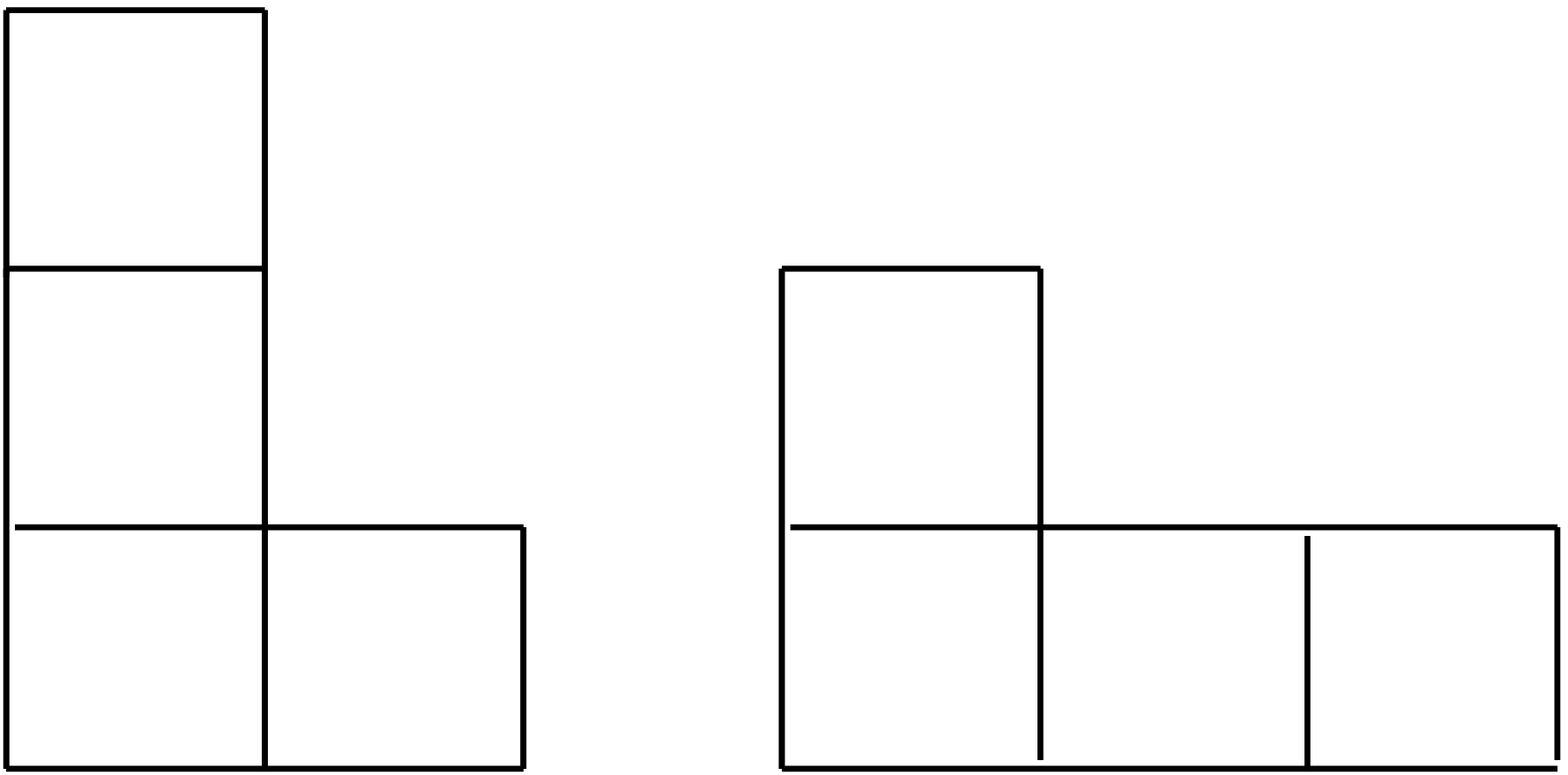}{100mm}{0}{The
Young diagrams $D_2$ and $D_3$} Here we have:
$p_{3,4}=18t^4+72t^5+162t^6+\mathcal{O}(t^7)$. The first space with
dimension 18 is the space of Lie brackets $[x_i,[x_j,[x_k,x_l]]]$ on
3 letters $x_1,x_2,x_3$. As $\gl_3$-module, it is the direct sum of
two representations with the Young diagrams $D_2$ and $D_3$ (see
Figure~3). The $\gl_3$-module with the left diagram, $D_2$, has
dimension 3, and the representation with the right diagram, $D_3$,
has dimension 15. The coinduced $W_3$-module from the direct sum of
these two representations on the level 4 has the character
$18\sum_{k\ge 0}\frac{(k+1)(k+2)}{2}t^{k+4}$. Differ this character
from $p_{3,4}$. The difference is $18t^5+54t^6+\mathcal{O}(t^7)$
which again has the character of the coinduced module from
$D_2\oplus D_3$ but on the level 5.

We can make this computation for any $A_{n,k}^\ell$ from the
tableaux  (\ref{eq1_3}). These computations motivate the Conjecture
above.
\endcomment
Consider the Lie algebra $grA_n$. Consider its first component
$A_n/[A_n,A_n]$. Consider the image under the canonical projection
$p\colon A_n\to A_n/[A_n,A_n]$ of the kernel
$K_n=A_n[A_n,[A_n,A_n]]A_n$. Denote this image by $\mathcal{Z}$.
\begin{lemma}
$\mathcal{Z}$ belongs to the center of the Lie algebra $grA_n$.
\begin{proof}
It follows from Lemma 2.2.1.
\end{proof}
\end{lemma}
Denote the quotient Lie algebra $grA_n/\mathcal{Z}$ by
$\wtilde{gr}A_n$. It is $\mathbb{Z}_{\ge1}$-graded Lie algebra.
\begin{theorem*}
On each graded component of $\wtilde{gr}A_n$ there acts the Lie
algebra $W_n$ of polynomial vector fields on an $n$-dimensional
vector space. The bracket is $W_n$-equivariant.
\begin{proof}
Consider the Lie algebra $Der(A_n)$ of the derivations of $A_n$ {\it
as of associative algebra}. This Lie algebra can be easily
described: a derivation of a free associative algebra is uniquely
defined by its values on the generators $\{x_i\}_{i=1,\dots,n}$, and
these values may be arbitrary. This Lie algebra acts on $A_n$, and
it induces the action on any quotient like $A_n/[A_n,A_n]$,
$A_n/A_n[A_n,A_n]A_n$, etc. In particular, it acts on
$A_n/A_n[A_n,[A_n,A_n]]A_n$. The last space is isomorphic to
$\Omega^{even}(\mathbb{C}^n)_*$ as algebra (considered as an algebra
in degree 0) by Lemmas 2.1.1 and 2.1.2(iii). Then we have a map
$\varphi\colon Der(A_n)\to Der(\Omega^{even}(\mathbb{C}^n)_*)$.
Denote by $\Im$ the image of this map.
\begin{lemma}
The Lie algebra $\Im$ acts on $\wtilde{gr}A_n$, and the bracket in
$\wtilde{gr}A_n$ is $\Im$-equivariant.
\begin{proof}
We need to prove that $\Im$ acts on each consecutive quotient. For
this we need to prove that if we apply a derivation such that all
generators $x_i$ are mapped to $A_n[A_n,[A_n,A_n]A_n$ to a
$k$-commutator in $A_n$, the image will belong to
$(k+1)$-commutator. It follows (for any $k$) immediately from Lemma
2.2.1.
\end{proof}
\end{lemma}
Now we just should construct a Lie subalgebra isomorphic to $W_n$ in
$\Im$. This is the Lie algebra $W_n$ which acts in the natural way
on all even forms in $\Omega^{even}(\mathbb{C}^n)_*$. It is clear
that this subalgebra belongs to the image $\Im$. Indeed, each
derivation of $A_n/A_n[A_n,[A_n,A_n]]A_n$ is defined by its values
on the generators $\{x_i\}$. These values are defined up to
$A_n[A_n,[A_n,A_n]]A_n$. Take an arbitrary lift of each value in
$A_n$, we get a derivation of $A_n$. This speculation proves also
that the map $Der(A_n)\to Der(A_n/A_n[A_n,[A_n,A_n]]A_n)$ is
surjective, and $\Im=Der(\Omega^{even}(\mathbb{C}^n)_*)$. Consider
the Lie algebra $W_n$ acting on {\it commutative} forms
$\Omega^{even}(\mathbb{C}^n)$ in the natural way. Then it acts on
the quantized algebra $\Omega^{even}(\mathbb{C}^n)_*$ as well. In
particular, the canonical $W_n$ acts on $\tilde{gr}A_n$.
\end{proof}
\end{theorem*}

\subsection{Examples}
\begin{example}
The first grading component of the Lie algebra
$\wtilde{gr}^1(A_n)=\Omega^{even}/{\rm Im}d$. The bracket
$\Lambda^2(\wtilde{gr}^1(A_n))\to \wtilde{gr}^2(A_n)=[A_n,A_n]/
[A_n,[A_n,A_n]]=\Omega^{even+}_{closed}$ is the map
$\omega_1\bigwedge\omega_2\to (d\omega_1)\wedge (d\omega_2)$ which
is clearly $W_n$-equivariant.
\end{example}
\begin{example}
The following computation was made using MAGMA by Eric Rains.

Let $F_1=A_n$, and $F_k=[A_n,F_{k-1}]$ for $k>1$. Denote
$A_{n,k}=F_k/F_{k+1}$, and denote by $A_{n,k}^\ell$ the graded
component of $A_{n,k}$ consisting from the monomials of degree
$\ell$. Consider the bigraded Hilbert series for $A_n$:
\begin{equation}\label{eq1_2}
H_{n}=\sum_{\ell\ge 0,\ k\ge 1} \dim {A_{n,k}^\ell}u^k t^\ell
\end{equation}
For $n=2$ and for $n=3$ the bigraded Hilbert series are:
\begin{equation}\label{eq1_3}
\begin{aligned}
\ & H_2(u,t)=\\
&\ \ \ (  u)\\
& +( 2u)t\\
& +( 3u+ u^2)t^2\\
& +( 4u+2u^2+ 2u^3)t^3\\
& +( 6u+3u^2+ 4u^3+ 3u^4)t^4\\
& +( 8u+4u^2+ 6u^3+ 8u^4+ 6u^5)t^5\\
& +(14u+5u^2+ 8u^3+13u^4+15u^5+ 9u^6)t^6\\
&+(20u+6u^2+10u^3+18u^4+26u^5+30u^6+ 18u^7)t^7\\
&+(36u+7u^2+12u^3+23u^4+37u^5+54u^6+ 57u^7+ 30u^8)t^8\\
&+(60u+8u^2+14u^3+28u^4+48u^5+80u^6+108u^7+110u^8+56u^9)t^9\\
&+\mathcal{O}(t^{10})\\
& \\
&\\
 &H_3(u,t)=\\
 &\ \ \  (   u)\\
&+(  3u)t\\
& +(  6u+ 3u^2)t^2\\
& +( 11u+ 8u^2+ 8u^3)t^3\\
& +( 24u+15u^2+24u^3+ 18u^4)t^4\\
& +( 51u+24u^2+48u^3+ 72u^4+ 48u^5)t^5\\
&+(130u+35u^2+80u^3+162u^4+206u^5+116u^6)t^6\\
&+\mathcal{O}(t^7)
\end{aligned}
\end{equation}

We will find the consecutive quotients as $W_n$-modules of the class
$\mathcal{C}$, for small $k$, and $n=2,3$. We have already proved in
the paper that $A_{n,2}$ is a $W_n$-module. Consider the space
$A_{n,3}$ for $n=2$. We are going to show that the dimensions $\dim
A_{2,3}^{\ell}$ are exactly like the character of a $W_2$-module.
Indeed, we know from (\ref{eq1_3}) that $\dim
A_{2,3}^\ell=2(\ell-2)$ for $3\le \ell\le 9$. Consider a
$W_2$-module coinduced from a Young diagram $D_1$ showed in Figure 2
with a 2-dimensional $\gl_2$-module on the level $\ell=3$. Then on a
level $\ell$ this $W_2$-module should have dimension $2
S_{\ell-3,2}$ where $S_{k,2}$ is dimension of symmetric polynomials
in 2 variables of degree $k$. We have: $S_{k,2}=k+1$. This shows
that the spaces $A_{2,3}^\ell$ have the character of a $W_2$-module
for $\ell\le 9$.

Consider now $A_{2,4}$. Here the polynomial for $A_{2,4}^\ell$ is
$3t^4+8t^5+13t^6+18t^7+23t^8+28t^9+\mathcal{O}(t^{10})$. First take
the difference with $3t^4+6t^5+\dots+3(k+1)t^{k+4}+\dots$ which is
the character of a $W_2$-module. The difference is $\sum_{k\ge
0}2(k+1)t^{k+5}$ which is a character of $W_2$-module.

The case $A_{2,5}$ is analogous.

It is more interesting to consider the case $n=3$. Consider
$A_{3,3}$. \sevafigc{pasha2.eps}{40mm}{0}{The Young diagram $D_1$}
The polynomial for $A_{3,3}^\ell$ is
$p_{3,3}=8t^3+24t^4+48t^5+80t^6+\mathcal{O}(t^7)$. The first
coefficient $8$ is the space generated by Lie words
$[x_i,[x_j,x_k]]$ on 3 letters $x_1,x_2,x_3$ (some of $i,j,k$ may
coincide). This space has dimension 8, and as $\gl_3$-module it is
corresponded to the Young diagram $D_1$ (see Figure~2). The
irreducible $\gl_3$-module corresponded to this Young diagram has
dimension 8, and it can be easily computed by the character formula.
Consider the $W_3$-module coinduced from this $\gl_3$-module. This
coinduced character is $8\sum_{k\ge 0}\frac{(k+1)(k+2)}{2}t^{k+3}$.
It is exactly our $p_{3,3}$ up to $t^6$.

Finally, consider $A_{3,4}$. \sevafigc{pasha3.eps}{100mm}{0}{The
Young diagrams $D_2$ (left) and $D_3$ (right)} Here we have:
$p_{3,4}=18t^4+72t^5+162t^6+\mathcal{O}(t^7)$. The first space with
dimension 18 is the space of Lie brackets $[x_i,[x_j,[x_k,x_l]]]$ on
3 letters $x_1,x_2,x_3$. As $\gl_3$-module, it is the direct sum of
two representations with the Young diagrams $D_2$ and $D_3$ (see
Figure~3). The $\gl_3$-module with the left diagram, $D_2$, has
dimension 3, and the representation with the right diagram, $D_3$,
has dimension 15. The coinduced $W_3$-module from the direct sum of
these two representations on the level 4 has the character
$18\sum_{k\ge 0}\frac{(k+1)(k+2)}{2}t^{k+4}$. Subtract this
character from $p_{3,4}$. The difference is
$18t^5+54t^6+\mathcal{O}(t^7)$ which again has the character of the
coinduced module from $D_4\oplus D_5$ on the level 5 (see Figure~4).
The irreducible $\gl_3$-module corresponding to the Young diagram
$D_4$ has dimension 3, and the irreducible $\gl_3$-module
corresponding to $D_5$ has dimension 15. The number of boxes in the
Young diagram should be equal to the length of the words in the
corresponding representation of $\gl_3$, that is, to the level
$\ell$. \sevafigc{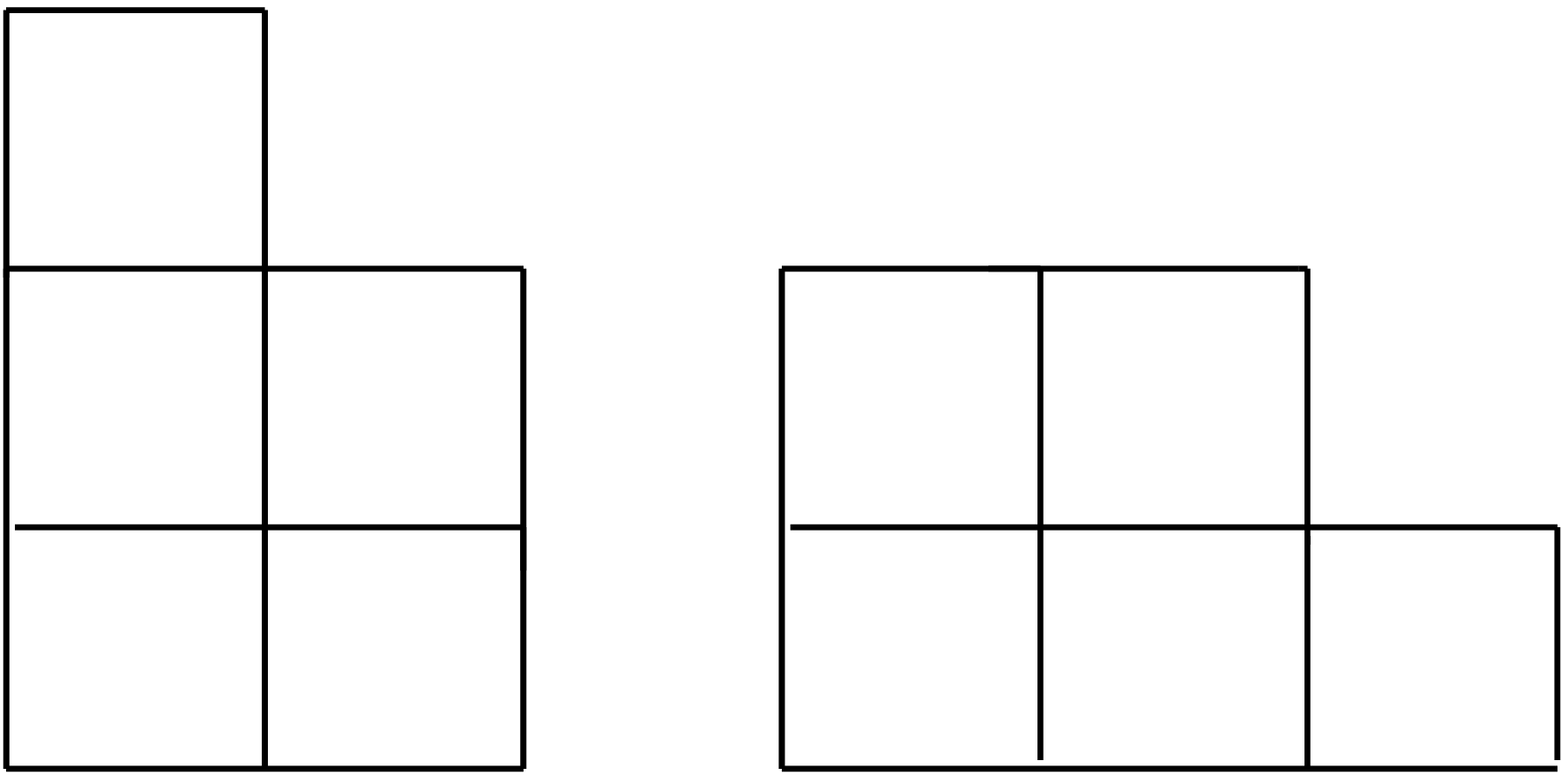}{100mm}{0}{The Young diagrams $D_4$
(left) and $D_5$ (right)}
\end{example}

\comment

\section{The quotient
$[A_n,[A_n,A_n]]/[A_n,[A_n,[A_n,A_n]]]$}
\subsection{ A generalization of Lemma 1.2.1}
\subsubsection{}
Let $A$ be an associative algebra. Consider the map
\begin{equation}\label{eq5_1}
\varphi\colon A\otimes [A,A]\to [A,[A,A]]
\end{equation}
which maps $a\otimes[b,c]\mapsto [a,[b,c]]$. It is clear that this
map is surjective. What is its kernel?

First of all, it contains the "Jacobi identities":
\begin{equation}\label{eq5_2}
\varphi(a\otimes [b,c]+b\otimes [c,a]+c\otimes [a,b])=0
\end{equation}
It is also clear that
\begin{equation}\label{eq5_3}
\varphi([a,b]\otimes [c,d]+[c,d]\otimes[a,b])=0
\end{equation}
It is clear that these two equations do {\it not} exhaust the kernel
of the map $\varphi$.
\begin{lemma}
The cyclic sum
\begin{equation}\label{eq5_4}
[a_1,[a_2,a_3a_4]]+[a_2,[a_3,a_4a_1]]+[a_3,[a_4,a_1a_2]]+[a_4,[a_1,a_2a_3]]
\end{equation}
is 0. \qed
\end{lemma}
Thus we have the relation
\begin{equation}\label{eq5_5}
\varphi(a_1\otimes[a_2,a_3a_4]+a_2\otimes[a_3,a_4a_1]+a_3\otimes[a_4,a_1a_2]+a_4\otimes[a_1,a_2a_3])=0
\end{equation}
We formulate the following conjecture, which is the direct analog of
Lemma 1.2.1.
\begin{conjecture}
For $A=A_n$, we have the isomorphism
\begin{equation}\label{eq5_6}
\overline{\varphi}\colon (A\otimes [A,A])/({\it relations}\
(\ref{eq5_2}),(\ref{eq5_3}),(\ref{eq5_5}))\eqto [A,[A,A]]
\end{equation}
\end{conjecture}
We conjecture the following cohomological meaning of the kernel
$\overline{\varphi}$ for a general associative algebra $A$. Consider
the functor (generalizing the cyclic functor $A\mapsto A/[A,A]$)
which maps an associative algebra $A$ to the vector space
$[A,A]/[A,[A,A]]$. Denote this functor by $\theta_2(A)$ (we denote
$\theta_1(A)=A/[A,A]$). We now can define the derived functors of
$\theta_2$ as follows: replace the algebra $A$ by a free associative
dg algebra $R^\mb(A)$, which is a resolution of $A$ in the class of
free associative dg algebras. Then compute the cohomology
$H^i(\theta_2(R^\mb(A)))=HC^i_{(2)}(A)$. One can prove that these
spaces do not depend on the choice of the resolution $R^\mb(A)$. We
call them {\it (2)-cyclic homology} of $A$. We conjecture that the
kernel of the map $\overline{\varphi}$ from (\ref{eq5_6}) is the
first (2)-cyclic homology group $HC^1_{(2)}(A)$. In particular, it
is compatible with the conjecture that for a free $A$ the map
$\overline{\varphi}$ is an isomorphism. Actually, these two
conjectures are equivalent.
\subsubsection{}
It is clear that the map $\overline{\varphi}$ is a map of
$A$-modules where $A$ is considered as the Lie algebra. Then for
$A=A_n$, the isomorphism (\ref{eq5_6}) should induce the isomorphism
of $A$-coinvariants.
\begin{lemma}
The coinvariants module of $(A\otimes [A,A])/({\it relations}\
(\ref{eq5_2}),(\ref{eq5_3}),(\ref{eq5_5}))$ is isomorphic to
$(A/[A,A]\otimes [A,A]/[A,[A,A]])/({\it relations}\
(\ref{eq5_2}),(\ref{eq5_3}),(\ref{eq5_5}))$.
\begin{proof}
We have: $[x,a\otimes[b,c]]=[x,a]\otimes [b,c]+a\otimes[x,[b,c]]$.
Let us simplify the r.h.s.: $[x,a]\otimes [b,c]=-[b,c]\otimes[x,a]$
by the relation (\ref{eq5_3}). Next, by (\ref{eq5_2}),
$-[b,c]\otimes[x,a]=a\otimes[[b,c],x]+x\otimes[a,[b,c]]$. Finally,
we get $[x,a\otimes[b,c]]=x\otimes [a,[b,c]]$. It proves that the
elements in $A\otimes [A,[A,A]]$ belong to the image of the action
of $A$ on $A\otimes [A,A]$. Then it is clear that $[A,A]\otimes
[A,A]$ also belongs to the image of the $A$-action.
\end{proof}
\end{lemma}
Notice that we used only relations (\ref{eq5_2}) and (\ref{eq5_3})
in the proof and did not use the relation (\ref{eq5_5}).
\begin{corollary}
If Conjecture 5.1 is true, then for $A=A_n$ one has the isomorphism
\begin{equation}\label{eq5_7}
\overline{\varphi}\colon(A/[A,A]\otimes [A,A]/[A,[A,A]])/({\it
relations}\ (\ref{eq5_2}),(\ref{eq5_3}),(\ref{eq5_5}))\eqto
[A,[A,A]]/[A,[A,[A,A]]]
\end{equation}
\qed
\end{corollary}
\subsection{ The double forms $\Omega(d,D)$}
\subsubsection{}
\begin{lemma}
\begin{itemize}
\item[(i)]
Let $A^\mb$ be an associative dg algebra (the differential of degree
+1 obeys $d^2=0$ and satisfies the graded Leibniz rule). Define a
new product on $A^\mb$ (not graded) as follows:
\begin{equation}\label{eq6_1}
\omega_1\star\omega_2=\omega_1\cdot \omega_2+(-1)^{\deg
\omega_1}d\omega_1\cdot d\omega_2
\end{equation}
for homogeneous $\omega_1$ and $\omega_2$. Then this product is
associative,
\item[(ii)]
Let $A^{\mb}$ be a graded algebra with two differentials $d$ and $D$
of degree +1 which both satisfy $d^2=D^2=0$ and the graded Leibniz
rule, and which skew-commute: $dD+Dd=0$. Then the product
\begin{equation}\label{eq6_2}
\omega_1\star\omega_2=\omega_1\cdot\omega_2+(-1)^{\deg
\omega_1}d\omega_1\cdot d\omega_2+(-1)^{\deg \omega_1}D\omega_1\cdot
D\omega_2-(dD\omega_1)\cdot (dD\omega_2)
\end{equation}
(for homogeneous $\omega_1$ and $\omega_2$) is associative.
\end{itemize}
\begin{proof}
We show how to deduce (ii) from (i). The statement (i) can be
checked directly.

Consider the algebra $A^\mb_*$ with the product from the part (i) of
the Lemma, but now $A^\mb$ is endowed also with the differential
$D$. We used only $\mathbb{Z}_2$-grading in the statement (i). Then
we first prove that the pair $(A^\mb_*,D)$ satisfies the conditions
of the statement (i). We now denote the product from the statement
(i) by $\star_0$. We need to prove that
\begin{equation}\label{eq6_3}
D(\omega_1\star_0\omega_2)=D\omega_1\star_0\omega_2+(-1)^{\deg\omega_1}\omega_1\star_0
D\omega_2
\end{equation}
This is a direct check. We use that $dD+Dd=0$. Now we are in the
conditions (i) and we can write the associative product
$\omega_1\star\omega_2=\omega_1\star_0\omega_2+(-1)^{\deg\omega_1}D\omega_1\star_0
D\omega_2$. The latter coincides with the product in the item (ii)
of lemma because $\deg D\omega_1=\deg\omega_1+1$.
\end{proof}
\end{lemma}
\subsubsection{}
Now consider the map $\varphi_n(d,D)$ of the free associative
algebra $A_n=<x_1,\dots,x_n>$ to the algebra $\Omega^{even}(d,D)$.
We have: $\varphi_n(d,D)(x_i)=x_i$, and then we define the map
$\varphi_n(d,D)$ on all monomials of $A_n$ as a map of algebras.

In the case of one differential, the map
$\varphi_n=\varphi_n(d)\colonA_n\to\Omega^{even}$ is surjective (see
Lemma 3.1.1). For the map $\varphi_n(d,D)$ the analogous statement
is not true.

\bigskip

{\bf Question.} How to define the image of $\varphi_n(d,D)$ in
$\Omega^{even}(d,D)$?

Some very naive conjecture about this question could be the
following:

The first consecutive quotient $A_n/[A_n,A_n]$ is very big while the
next quotients are sufficiently small (see Section 4). Then we could
conjecture that all consecutive quotients of the commutant
filtration on $\Omega(d,D)_*^{even}$ starting from
$[\Omega(d,D)_*^{even},\Omega(d,D)_*^{even}]/[\Omega(d,D)_*^{even},[\Omega(d,D)_*^{even},\Omega(d,D)_*^{even}]]$
{\it coincide} with the corresponding quotients for the image of
$A_n$ in $\Omega(d,D)_*^{even}$. In particular, the commutator
$[\Omega(d,D)_*^{even},\Omega(d,D)_*^{even}]$ is the same as for the
image of $A_n$ in $\Omega(d,D)_*^{even}$.

Maybe this conjecture is not true in such a strong form, but this
phenomena starts from some $k$-th commutator. But I think that it
starts from the commutator
$[\Omega(d,D)_*^{even},\Omega(d,D)_*^{even}]$.

\subsubsection{}
Now let me formulate some conjecture on the kernel of the map
$\varphi_n(d,D)\colonA_n\to\Omega(d,D)_*^{even}$. Firstly let us
notice that for any $k\ge 0$ the $k$-th term of the commutator
filtration on $\Omega(d,D)_*^{even}$ is not 0. (In the case of 1
differential  we have in Lemma 3.1.2(ii) that the commutator
$[\Omega^{even}(d)_*,[\Omega^{even}(d)_*,\Omega^{even}(d)_*]]=0$).

Here we have several first commutators in $\Omega(d,D)_*^{even}$:
\begin{equation}\label{eq6.5}
\begin{aligned}
\ &[\omega_1,\omega_2]=2(d\omega_1\cdot d\omega_2+D\omega_1\cdot
D\omega_2)\\
&[\omega_3,[\omega_1,\omega_2]]=4(D\omega_3\cdot D(d\omega_1\cdot
d\omega_2)+d\omega_3\cdot d(D\omega_1\cdot D\omega_2))\\
&[\omega_4,[\omega_3,[\omega_1,\omega_2]]]=8(d\omega_4\cdot
dD\omega_3\cdot D(d\omega_1\cdot d\omega_2)+D\omega_4\cdot
Dd\omega_3\cdot d(D\omega_1\cdot D\omega_2))
\end{aligned}
\end{equation}
We see that the last line "has more symmetries" than the previous
two. We are going to write down a linear combination of the elements
of the form $[m_4,[m_3,[m_2,m_1]]]$, $m_i\inA_n$, which is mapped by
$\varphi(d,D)$ to 0.

\begin{lemma}
Let $m_1,m_2,m_3,m_4\in A_n$. Denote $R(m_1,m_2,m_3,m_4)=
[m_4,[m_3,[m_1,m_2]]]-[m_1,[m_3,[m_4,m_2]]]+[m_2,[m_3,[m_1,m_4]]]$.
Then the map $\varphi(d,D)\colonA_n\to\Omega(d,D)_*^{even}$ maps
$R=R(m_1,m_2,m_3,m_4)-R(m_1,m_3,m_2,m_4)$ to 0.
\begin{proof}
It is a direct check.
\end{proof}
\end{lemma}

Let us notice that there are 6 linearly independent commutators
between $m_1,m_2,m_3,m_4$, and our relation $R$ has 6 terms.

\begin{conjecture}
The kernel of the map $\varphi(d,D)$ is $A_n\cdot\{R\}\cdot A_n$.
\end{conjecture}

It would be very interesting to check is it true or not.

\subsubsection{}
And the last question: we think that
$[A_n,[A_n,A_n]]/[A_n,[A_n,[A_n,A_n]]]$ is isomorphic to
$[\Omega(d,D)_*^{even},[\Omega(d,D)_*^{even},\Omega(d,D)_*^{even}]]/[\Omega(d,D)_*^{even},[\Omega(d,D)_*^{even},[\Omega(d,D)_*^{even},\Omega(d,D)_*^{even}]]]$
under the isomorphism $\varphi(d,D)$ and Conjecture 5.2.2. This is
the role of the algebra $\Omega(d,D)_*^{even}$ in the theory.
\endcomment

\subsection*{Acknowledgements.}We are grateful to Pavel Etingof, Giovanni
Felder, Anton Khoroshkin,  and especially to Misha Movshev and to
Eric Rains for many discussions. B.F. is grateful to the grants RFBR
05-01-01007, SS 2044.2003.2, RFBR-JSPS 05-01-02934, and also to the
Russian Academy of Science project "Elementary particle physics and
fundamental nuclear physics" for partial financial support. B.Sh. is
thankful to the ETH (Zurich), and to the University of Luxembourg,
where some ideas of this paper were invented, for excellent working
conditions. He is also grateful to the ETH Research Comission grant
TH-6-03-01 and to the research grant R1F105L15 of the University of
Luxembourg for partial financial support.

\bigskip
\bigskip
B.F:\\
Landau Institute for Theoretical Physics, 1A prospekt Akademika
Semenova, Chernogolovka, Moscow Region, 142432 RUSSIA\\
and\\
Independent University of Moscow, 11 Bol'shoj Vlas'evskij per.,
119002 Moscow RUSSIA\\
{\it e-mail}: {\tt bfeigin@gmail.com}\\
\\
B.Sh:\\
Faculty of Science, Technology and Communication, Campus
Limpertsberg, University of Luxembourg,
162A avenue de la Faiencerie, L-1511 LUXEMBOURG\\
{\it e-mail}: {\tt borya$\_$port@yahoo.com}

\end{document}